\documentclass[]{scrartcl}
\pdfoutput=1
\usepackage{ae,lmodern}
\usepackage[english]{babel}
\usepackage[utf8]{inputenc}
\usepackage[T1]{fontenc}
\usepackage{mathtools}
\usepackage{amsfonts}
\usepackage{amssymb}
\usepackage{mathrsfs}
\usepackage{dsfont}
\allowdisplaybreaks

\usepackage[
backend=biber,
style=trad-alpha,
eprint=true,
doi=false,
url=false,
isbn=false,
]{biblatex}

\usepackage{geometry}
\usepackage{abstract}
\usepackage{microtype}

\usepackage[bookmarks=true, breaklinks]{hyperref}
\hypersetup{
	colorlinks=true,
	linkcolor=black,
}

\usepackage{amsthm}
\usepackage{footnotebackref}
\hypersetup{pdfauthor={Charles Bordenave, Joffrey Mathien},pdftitle={Cutoff for geodesic paths on hyperbolic manifolds}}

\newcommand{\DTV}{\mathrm{d_{TV}}}
\newcommand{\cB}{\mathcal B}
\newcommand{\cL}{\mathcal L}

\newcommand{\cD}{\mathcal D}

\newcommand{\dP}{\p}
\newcommand{\dH}{\mathds H}
\newcommand{\IND}{\mathbf{1}}

\newcommand{\R}{\mathds{R}}
\newcommand{\C}{\mathds{C}}
\newcommand{\Rp}{\R^{+}}
\newcommand\limit[2]{\overset{#1}{\underset{#2}{\longrightarrow}}}
\newcommand\limitninf[1]{\limit{#1}{n \to \infty}}

\newcommand{\p}{\mathbb{P}}
\newcommand{\E}{\mathbb{E}}
\newcommand{\eps}{\varepsilon}
\renewcommand{\d}{\mathrm{d}}

\DeclareMathOperator{\Supp}{Supp}
\DeclareMathOperator{\Vol}{Vol}
\DeclareMathOperator{\Diam}{Diam}

\newtheorem{theorem}{Theorem}[section]
\newtheorem{proposition}[theorem]{Proposition}

\newtheorem{corollary}[theorem]{Corollary}
\newtheorem{definition}[theorem]{Definition}

\title{Cutoff for geodesic paths 
on hyperbolic manifolds}
\author{Charles Bordenave\thanks{I2M, Aix-Marseille Université, CNRS. email:  \href{mailto:charles.bordenaven@univ-amu.fr}{\texttt{charles.bordenave@univ-amu.fr}}} \quad   \quad  Joffrey Mathien\thanks{I2M, Aix-Marseille Université, CNRS. email:  \href{mailto:joffrey.mathien@univ-amu.fr}{\texttt{joffrey.mathien@univ-amu.fr}}}}
\date{}

\numberwithin{equation}{section}
\numberwithin{figure}{section}

\begin{document}
	\maketitle
	
	\begin{abstract}
We establish new instances of the cutoff phenomenon for geodesic paths and for the Brownian motion on compact hyperbolic manifolds. We prove that for any fixed compact hyperbolic manifold, the geodesic path started on a spatially localized initial condition exhibits cutoff. Our work also extends results obtained by Golubev and Kamber on hyperbolic surfaces of large volume to any dimension. Our proof  builds upon a spectral strategy introduced by Lubetzky and Peres for Ramanujan graphs and on a detailed spectral analysis of the spherical mean operator. 
	\end{abstract}
	
	%\tableofcontents
	
	\section{Introduction}
	
		\paragraph{Cutoff Phenomenon} For an ergodic dynamical system, the cutoff describes an abrupt transition to equilibrium. Historically, it was introduced in seminal works of  Diaconis, Shahshahani and Aldous \cite{diaconis1981generating,aldous1983mixing,aldous1986shuffling} for card shuffling and other random walks on finite groups. 
	There are now numerous examples of Markov chains and Markov processes where cutoff has been established, see \cite{diaconis1996cutoff,MR3726904} but despite some recent efforts such as \cite{MR2375599,MR3650406,MR4476123,MR4780485}, a satisfactory theory of cutoff is however still missing. Also, most of current examples are on finite spaces with some exceptions such as \cite{MR2375599,MR3201989,Golubev,10.1214/20-EJP417,Boursier_2022}. In this paper, we study cutoff for classical processes on compact hyperbolic manifolds and we extend a spectral strategy which has been introduced in \cite{MR3558308} and further developed in \cite{Golubev,MR4167012,MR4476123,MR4400017}.

		The formal definition of cutoff that we shall use is as follows. Let $(X^\delta,\cB^\delta,\pi^\delta)$ be a sequence of probability spaces indexed by a parameter $\delta \in \mathcal \cD \subset \R$. Let $\delta_\star$ be an accumulation point of $\cD$. For example, for a card shuffling problem, we may take $\delta = n$ where $n$ is the number of cards in the deck, $\delta_\star = \infty$, $X^\delta$ is the symmetric group of order $n$  and $\pi^\delta$ is the uniform (or Haar) probability measure on $X^\delta$.  Next, for each $\delta \in \cD$, we consider a sequence of probability measures  $(\mu^\delta_t)_{t \geq 0}$, indexed by a real or integer $t \geq 0$. In the card shuffling example, for integer $t \geq 0$, $\mu^\delta_t$ would be the law of the deck of cards after $t$ shuffles. We assume that $(\mu^\delta_t)$ converges to equilibrium, that is for all $\delta \in \cD$, 
	\begin{equation}\label{eq:conveq}
		\lim_{t \to \infty} \mu_t^\delta = \pi^\delta,
		\end{equation}
	where the convergence is in total variation : for any probability measures $\mu,\nu$ on a measurable space $(X,\cB)$, their distance in total variation is
	\begin{equation}\label{def:TV}
	\DTV (\mu,\nu)  = \sup_{ B \in \cB} | \mu(B) - \nu(B) | =  \frac 1 2 \int \left| \frac{d \mu}{d \nu}(x) - 1    \right| d\nu(x)  \in [0,1].
	\end{equation}
We have $\DTV(\mu,\nu) = 1$ if and only $\mu$ and $\nu$ are mutually singular. To define cutoff, we consider an initial condition $\mu^\delta_0$ which is nearly singular with respect to $\pi^\delta$: 
$$
\lim_{\delta \to \delta_\star} \DTV( \mu^\delta_0 , \pi^\delta) = 1.
$$
In the card shuffling example, the initial condition will usually be deterministic : $\mu_0^\delta$ is a Dirac mass at a given permutation of cards.  Then, the cutoff occurs when the map $t \to \DTV (\mu^\delta_t,\pi^\delta)$ becomes close to a step function  as $\delta \to \delta_\star$. 	
	\begin{definition}{(Cutoff)}\label{defi:cutoff}
 The sequence of measures $\mu_t^\delta$ exhibits cutoff  if there exists a sequence $(t_\star^\delta)_{\delta \in \cD}$  in $\Rp$ such that for any $\eps \in (0,1)$, 
	$$
	\lim_{\delta \to \delta_\star}  \inf_{t \leq (1-\eps) t_\star^\delta} \DTV (\mu^\delta_t,\pi^\delta) = 1 \quad \hbox{ and } \quad  	\lim_{\delta \to \delta_\star}  \sup_{t \geq (1 + \eps) t_\star^\delta} \DTV (\mu^\delta_t,\pi^\delta) = 0.
	$$
	We then say that the cutoff occurs at time $t_\star^\delta$. 
	\end{definition}
 For Markov chains, this definition of cutoff turns out to be equivalent to the usual definition of cutoff with mixing times (which we have skipped here), see \cite{MR3726904}. However, we retain the name \textit{mixing time} for $t_\star^\delta$. The choice of the total variation distance is the most natural from a probabilistic perspective and our definition is tuned for this distance. It is also possible to define cutoff for other distances as in \cite{MR2375599,MR3558308}. 
 
 We shall write $\cL(Z)$ for the law of a random variable $Z$. With a little abuse, we shall say that $Z_t^\delta$ exhibits cutoff if $\cL(Z_t^\delta)$ does. %In the same vein, we may write $\DTV( Z_1,Z_2)$ in place of $\DTV(\cL( Z_1),\cL(Z_2))$.

\paragraph{Geodesic process and spherical means on hyperbolic manifolds}
Let $M$ be a compact connected {oriented} manifold without boundary of constant sectional curvature $-1$ -- hereafter simply called {\em compact  hyperbolic manifold} --  of dimension $d \geq 2$.  We denote by  $T_x M$ the tangent space at $x \in M$, $T^1_x\hspace{-1pt}M$ the unit $(d-1)$-sphere of $T_x M$ and by $TM$, $T^1\hspace{-2pt}M$ the tangent and unit tangent bundles over $M$. The volume form associated with the Riemannian metric on $M$ is denoted $\Vol$ and $\pi = \Vol / \Vol(M)$ is the probability measure associated with $\Vol$. We denote by $\sigma_x$ the uniform probability measure on $T^1_x\hspace{-1pt}M$ so that the Liouville probability measure on $T^1\hspace{-2pt}M$ can be written as $d L (z)  = d\pi(x) d\sigma_x(v)$, with $z = (x,v) \in  T^1\hspace{-2pt}M$.

In this paper, we consider the geodesic flow $\gamma(z, t)$ for $z \in T^1\hspace{-2pt}M, t \in \Rp$,  and our goal is to show that, from a probabilistic point of view, it exhibits some cutoff phenomenon. More precisely, we study the following {\em geodesic process}. Pick a point $Z_0 = (X_0,v_0) \in T^1\hspace{-2pt}M$ at random where $X_0 \in M$ has law $\mu_0$, and, given $X_0$,  $v_0 \in T^1_{X_0}M$ has conditional law $\sigma_{X_0}$. Then define \[Z_t = (X_t,v_t) = \gamma (Z_0,t), \quad \mu_t = \mathcal{L}(X_t).\] 

Note that our geodesic process it is not what is generally called the geodesic random walk, studied for example in \cite{MR758647,Golubev}. Here, the randomness of the process $Z_t$ is only through the initial condition $Z_0$.

The process $X_t$ is associated with a family of bounded operators $A_t$ acting on $L^2(M)$: 
\[A_tf(x) = \int_{T^1_xM} f(\gamma((x,v),t)) d\sigma_x(v).\]
The operator $A_t$ is called the {\em spherical mean of radius $t$}. If $f_0$ is the density of $\mu_0$ with respect to $\pi$ and $f_0 \in L^2(M)$, then $\mu_t$ is also absolutely continuous with respect to $\pi$, and its density is given by $A_tf_0$.

Geometrically, our geodesic process can also be described as follows. The universal cover of $M$ is the $d$-dimensional hyperbolic space $\mathds{H}_d$ which is naturally endowed with the hyperbolic measure $\ell_{d}$. 
When $(M,X_0)$ is lifted to $(\mathds{H}_d,\widetilde X_0)$, the lift of $X_t$, say $\widetilde X_t$, is sampled according to the projection  of $\ell_d$ to the sphere in $\mathds{H}_d$ of radius $t$ and center $\widetilde X_0$. 

It is a general fact that $\mu_t$ converges to $\pi$ as $t \to \infty$, see \cite{MR200874,MR1306567,pjm}. In this work, we choose $\mu_0$ to be uniform on a small ball around a fixed point $x \in M$. This choice is motivated by the fact that on the one hand we need a measure which is continuous with respect to $\pi$ to have convergence in total variation (which would not have been the case with a Dirac mass). On the other hand, this choice still describes well what happens around $x$. We prove that $\mu_t$ exhibits cutoff when $ \mu_t = \mu_t^\delta$ is parametrized according to two settings which we describe below. Finally, we observe that whereas $\mu_t$ converges in total variation, this is not the case for $\cL(Z_t)$ toward the Liouville measure $L$. In fact, due to the $L$-invariance and invertibility of the geodesic flow, we have $\DTV( \cL(Z_t) , L) = \DTV ( \cL(Z_0),L)$. We refer to \cite{MR896798,MR3338007} and references therein for studies of the mixing rate of the geodesic flow.
	
\paragraph{Cutoff with a peaked initial condition}
For $\delta >0$ and $x \in M$, let $\mu_0^\delta$ be the uniform probability measure on the ball of radius $\delta$ around $x$. As above, let $X_t^\delta$ be the geodesic process with initial probability measure $\mu_0^\delta$.

\begin{theorem}{}\label{th:D1}
	Let $M$ be a compact connected hyperbolic manifold of dimension $d$. For all $x \in M$, as $\delta \to 0$, the geodesic process $X^\delta_t$ exhibits  cutoff at time $t^\delta_\star = -\frac{\ln \delta}{d-1}$.
\end{theorem}

In plain words, this result asserts that the position of the geodesic flow on $T^1\hspace{-2pt}M$ starting from a peaked initial condition on $M$ and uniform on the tangent sphere satisfies the cutoff at a time which depends only on the dimension of $M$. The proof will also show that the cutoff is uniform over the initial condition $x \in M$ and that the cutoff window is order $1$. That is for any $\eta \in (0,1)$, there exists $T  = T(\eta)>0$ such that for all $x \in M$,
$$ \inf_{t \leq  t_\star^\delta - T} \DTV (\mu^\delta_t,\pi)  \geq 1 - \eta \quad \hbox{ and } \quad  	 \sup_{t \geq  t_\star^\delta  + T} \DTV (\mu^\delta_t,\pi ) \leq \eta.
	$$

From a probabilistic perspective, one reason why Theorem \ref{th:D1} is of interest is that it is an example of cutoff where the underlying probability space $(M,\pi)$ does not depend on $\delta$. The only other example known to the authors is \cite{10.1214/20-EJP417} for the Ornstein-Uhlenbeck process driven by a slowed-down Brownian motion (or Lévy process). Here, the mixing mechanism seems very different.

\paragraph{Laplace-Beltrami operator} Before introducing our next results, we recall some basic facts on $\Delta = - \mathrm{div} (\nabla)$, (minus) the Laplace-Beltrami operator on the manifold $M$. We refer to \cite{Stri,Buser, Chavel}. 	

	 On a compact connected hyperbolic manifold (without boundary), $\Delta$ extends to a positive self-adjoint operator with  discrete spectrum $\mathfrak S(\Delta)$ with countably many eigenvalues which we denote,  counting multiplicities, 
	 $$0 = \lambda_0(M) <\lambda_1 (M) \leq \lambda_2 (M) \leq  \dots,$$
see Subsection \ref{subsec:UB} below. 
On $\dH_d$,  the Laplace-Beltrami operator has continuous $L^2$-spectrum on the half-line
$ \left[ \sigma_d  , \infty \right)$ with 
\begin{equation}\label{eq:defsigmad}
\sigma_d =  \left(\frac{d-1}{2} \right)^2,
\end{equation} see \cite{10.4310/jdg/1214429509}.
	
 If $M_n$ is a sequence of compact hyperbolic manifolds of common dimension $d$ and with $\Vol(M_n) \to \infty$, it is known that 
	\begin{equation}\label{eq:AB}
\limsup_n	\lambda_1 (M_n) \leq \sigma_d,
	\end{equation}
see \cite[Corollary 2.3]{Cheng1975} (or \cite{MR365408} for the special case of surface).
	In the same spirit as \cite{MR3558308,Golubev,MR4167012,MR4476123,MR4400017}, we shall now prove that sequence of hyperbolic manifolds of increasing volume and with few Laplacian eigenvalues in the interval $[0, \sigma_d)$, the geodesic process exhibits cutoff.	In view of \eqref{eq:AB} such manifolds have nearly-optimal spectral edge behavior.

\paragraph{Cutoff for hyperbolic manifolds with a near optimal spectral gap}

We consider sequences $(M_n)_{n\geq 0}$ of compact connected hyperbolic manifolds without boundary all of fixed dimension $d \geq 2$. We denote by $\pi_n$ the normalized volume measure. We fix $\delta >0$. For $ x \in M_n$, we let $\mu^n_{x,0}$ be the uniform probability measure on the ball of radius $\delta$ and radius $x$. We denote by $Z^n_{x,t} = (X^n_{x,t},v^n_{x,t})$ the geodesic process with initial position $X^n_{x,0}$ sampled according to $\mu^n_{x,0}$. Finally, let $\mu^n_{x,t}$ be the law of $X^n_{x,t}$.  

We set $V_n = \Vol(M_n)$. In particular if $d = 2$, we have $V_n = 4 \pi (g-1)$ where $g=g_n$ is the genus of the surface $M_n$. We denote by $\rho(x)$ the injectivity radius at $x \in M_n$. 
\begin{theorem}{}\label{th:S1}
	Let $(M_n)$ be a sequence of compact hyperbolic $d$-dimensional manifolds such that $\lim_n V_n = \infty$ and  $\liminf_n \lambda_1(M_n) \geq \sigma_d$. Let $x_n \in M_n$ be a sequence of points such that $\liminf_n \frac{ \ln  \rho(x_n) }{\ln V_n}  = 0$. 
	Then the process $X^n_{x_n,t}$ exhibits cutoff at time $t^n_\star = \frac{\ln V_n}{(d-1)}$ when $n \to \infty$.
\end{theorem}

The analog of Theorem \ref{th:S1} for the Brownian motion and the geodesic random walk is proven in \cite{Golubev} for compact hyperbolic surfaces (that is $d=2$). The proof of Theorem \ref{th:S1} will be proven with a similar strategy. To our knowledge, there is currently no example of applications of Theorem \ref{th:S1} in dimension $d \geq 3$. In dimension $d = 2$, it was recently proven in \cite{MR4635304} that there exists a sequence of compact hyperbolic surfaces $M_n$ of increasing volume such that $\lim \lambda_1(M_n) =\sigma_2$. This is actually a generic assumption: it was announced in \cite{anantharaman2024spectralgaprandomhyperbolic} that a random hyperbolic surface of genus $g$ (and volume $V_n = 4\pi(g-1)$) sampled according to the Weil-Petersson measure on the moduli space of Riemannian metrics satisfies $\lim \lambda_1(M_n) =\sigma_2$ in probability as $n \to \infty$. More recently, it has been showed in \cite{MageePuderVanHandel} that this is also the typical behavior of covering spaces of closed hyperbolic surfaces.

\paragraph{Cutoff for hyperbolic manifolds with a Sarnak-Xue spectral density property}
The hypothesis $\liminf_n \lambda_1(M_n) \geq \sigma_d$ is very strong. However, it is possible to relax this spectral assumption as observed in related contexts in \cite{Golubev,MR4476123,MR4400017}. This relaxed spectral assumption is called the {\em  Sarnak-Xue spectral density property} in  \cite{MR4400017} after \cite{MR1131400}. More precisely, let $(M_n)_n$ be a sequence of hyperbolic $d$-dimensional manifolds as above with $\lim_n V_n = \infty$. For any $s \in [0,1)$, set \[I(s) = \limsup_n \frac{\ln \# \left\lbrace k :    \lambda_k(M_n) \leq \sigma_d (1-s^2)\right\rbrace}{ \ln V_n}.\]
We say that $(M_n)_n$ satisfies the \textit{Sarnak-Xue density property} if for all $s \in (0,1)$, $I(s) \leq 1-s$.

In dimension $d=2$, there are a few examples of hyperbolic surfaces satisfying the Sarnak-Xue property. For arithmetic modular surfaces $M_q = \Gamma(q) \backslash \dH^2$ , \cite[Corollary 1]{MR1131400} implies that the Sarnak-Xue property holds as $q \to \infty$. For a random hyperbolic surface sampled according to the Weil-Petersson measure, \cite[Theorem 2]{MR4408435} implies that the Sarnak-Xue density property holds with high probability as the genus $g$ goes to infinity (this could also be extracted from the estimates in \cite{MR4718401,anantharaman2023friedmanramanujanfunctionsrandomhyperbolic}). Similarly, for random $n$-covers of a hyperbolic surface, \cite[Theorem 1.8]{MR4431124} implies that the Sarnak-Xue density property holds with high probability. In all these examples, this property is proved as a consequence of Selberg trace formula and a sharp upper bound on the number of closed geodesic paths of length $L \leq A \ln V_n$ with $A \geq 2$.

\begin{theorem}{}\label{th:S2}
	Let $(M_n)$ be a sequence of hyperbolic $d$-dimensional manifolds such that  $\lim_n V_n = \infty$, $\liminf_n \lambda_1(M_n) >0$ and the Sarnak-Xue density property holds. Suppose further that there exists a measurable set $T_n \subset M_n$ such that $\lim_n \inf_{x \in T_n} \frac{\ln \rho(x)}{ \ln V_n} = 0$ and $\lim_n \pi_n(T_n) = 1$. Then there exists a measurable subset $S_n $ such that $\lim_n \pi_n(S_n) =1$ and for any sequence $x_n \in S_n$, the process $X^n_{x_n,t}$ exhibits cutoff at time $t^n_\star = \frac{\ln V_n}{d-1}$ when $n \to \infty$.
\end{theorem}

We note that a more quantitative upper bound on $1 - \pi_n(S_n)$ could be extracted from the proof of Theorem \ref{th:S2} Also, it is not complicated to combine the proofs of Theorems \ref{th:S1} and \ref{th:D1} to get a statement when $n \to \infty$ and $\delta \to 0$ simultaneously. For the sake of clarity and readability, we will not do so here.

\paragraph{Application to the Brownian motion}	
The Brownian motion $(W_t)_{t \geq 0}$ on a manifold $M$ can be defined as the process with infinitesimal generator $\frac{1}{2}\Delta$, where $\Delta$ is the Laplace-Beltrami operator.	There are analogs of Theorem \ref{th:S1} and Theorem \ref{th:S2} for the Brownian motions on a sequence of hyperbolic manifolds $(M_n)$. For $x \in M_n$, we denote below by $W^n_{x,t}$ the Brownian motion on $M_n$ started at $x$ at time $t =0$. The statements are read as follows.

\begin{corollary}{}\label{coro:Brownianmixing1}
Let $(M_n)$ be a sequence of compact hyperbolic $d$-dimensional manifolds such that $\lim_n V_n = \infty$ and  $\liminf_n \lambda_1(M_n) \geq \sigma_d$. Let $x_n \in M_n$ be a sequence of points such that $\liminf_n \frac {\ln  \rho(x_n)}{\ln V_n} = 0$. 
	Then the process $W^n_{x_n,t}$ exhibits cutoff at time $2  t^n_\star / (d-1)  = 2 \ln V_n$ when $n \to \infty$.
\end{corollary}

\begin{corollary}{}\label{coro:Brownianmixing2}
	Let $(M_n)$ a sequence of hyperbolic $d$-dimensional manifolds such that  $\lim_n V_n = \infty$, $\liminf_n \lambda_1(M_n) >0$ and the Sarnak-Xue density property holds. Suppose further that there exists a measurable set $T_n \subset M_n$ such that $\lim_n \inf_{x \in T_n} \frac{\ln \rho(x)}{\ln V_n} = 0$ and $\lim_n \pi_n(T_n) = 1$. Then there exists a measurable subset $S_n \subset M_n$ such that $\lim_n \pi_n(S_n) =1$ and for any sequence $x_n \in S_n$, the process $W^n_{x_n,t}$ exhibits cutoff at time $2  t^n_\star / (d-1) = 2 \ln V_n$ when $n \to \infty$.
\end{corollary}

For hyperbolic surfaces, Corollary \ref{coro:Brownianmixing1} was proven in \cite{Golubev}. For $n$-covers of a hyperbolic surface, Corollary \ref{coro:Brownianmixing2} is also proven in \cite{Golubev}.

\paragraph{Idea of the proof}
	To prove Theorem \ref{th:D1}-\ref{th:S2}, we aim to find more tractable quantities that give us good estimates of the total variation distance.
	Since the geodesic flow is of unit speed, it is easy to obtain a crude upper bound on the volume of the support of the process $X_n$, which leads to a lower bound on $t_\star$. Note that this bound does not require any assumption on the manifolds, which tells us in particular that the mixing time we find is minimal.
	On the other hand, we prove a kind of variance contraction for the density $f_t$ of $\mathcal L(X_t)$ to derive an upper bound on the mixing time. More precisely, we use the spectral decomposition of the operator $A_t$, obtained from that of the Laplace-Beltrami operator, to show an exponential decay of the total variation distance between $\mu_t$ and $\pi$.
	
	For Brownian motion (Corollaries  \ref{coro:Brownianmixing1} and \ref{coro:Brownianmixing2}),  the general strategy is quite similar.
	We can no longer bound the support of the process, but we can still use the drift of the hyperbolic Brownian motion to know that it is concentrated in a precise region of small volume with respect to $V_n$.
	For the upper bound, we note that the rotational invariance of the law of the Brownian motion ensures that, conditionally on the distance $D_t$ between the starting point and $W_t$, its law is the same as the geodesic process at time $D_t$. This allows us to refer back to the previous results.

\paragraph{Some perspectives} Theorem \ref{th:D1} gives an example of the cutoff phenomenon for a mixing dynamical system on a compact space started from an initial condition whose distance to equilibrium is large (here a spatially localized initial condition). It would be interesting to investigate the universality of this phenomenon for a large class of mixing dynamical systems. Other examples include baker's maps or Arnold's cat map. Geodesic paths on a manifold of non-constant non-negative curvature or chaotic billiards would be natural subjects of investigation. 

\paragraph{Organization of the rest of the paper}

In Section \ref{GF}, we study the geodesic process and prove our main technical results. We conclude this section by the proofs of Theorem \ref{th:D1}, Theorem \ref{th:S1} and Theorem \ref{th:S2}. In Section \ref{BM}, we study the Brownian motion and prove Corollary \ref{coro:Brownianmixing1} and Corollary \ref{coro:Brownianmixing2}.

	\paragraph{Acknowledgment} This work could not have been completed without many enlightening discussions with our friends and colleagues Adrien Boulanger and Jean Raimbault. We also thank Elon Lindenstrauss for suggesting to use H{\"o}rmander Theorem to improve Theorem \ref{th:D1} from "almost all $x \in M$" to "all $x \in M$". We thank Baptiste Dugué for pointing out an error in the value of $c_d$, which was off by a factor of one half, in the published version of Proposition \ref{prop:SpAt0}.
		
\section{Study of the geodesic flow} \label{GF}
This section is devoted to the proofs of Theorems \ref{th:D1}, \ref{th:S1}, \ref{th:S2}, which consist in studying the behavior of the geodesic flow and its mixing properties from a probabilistic point of view. 
The proofs of these three theorems are quite similar, and so we will summarize the common elements, before going into the computational details of each one separately. Our reasoning is based on two different arguments, one for the lower bound of the mixing time and one for the upper bound. Only the latter requires the spectral estimates assumed in Theorems \ref{th:S1}, \ref{th:S2}.

As explained earlier, we recall that the process is associated with a family of bounded operators $A_t$ acting on $L^2(M)$: \[A_tf(x) = \int_{T^1_xM} f(\gamma((x,v),t)) d\sigma_x(v).\] When necessary, we will use the exponent $A_t^M$ to make clear on which surface we are working.

In the following, we try to understand the process by studying $A_t$. For $x$ in a manifold $M$ and $\rho >0$, we write $B(x, \rho)$ for the ball of radius $\rho$ centered at $x$, and $S(x, \rho)$ for the corresponding sphere.

\subsection{Lower bound} In this section, we provide the tools to obtain the lower bound on the mixing time, that is the following proposition:
\begin{proposition}{}\label{prop:LB}
	We have that
	\begin{itemize}
		\item under the hypothesis of Theorem \ref{th:D1},
		\[\lim_{\delta \to 0}  \inf_{t \leq (1-\eps) t_\star^\delta} \DTV (\mu^\delta_t,\pi) = 1.\]
		\item under the hypothesis of Theorem \ref{th:S1} or Theorem \ref{th:S2},
		\[\lim_{n \to \infty}  \inf_{t \leq (1-\eps) t_\star^n} \DTV (\mu^n_t,\pi^n) = 1.\]
	\end{itemize}
\end{proposition}
The argument required is quite elementary, since it comes from a volume argument and does not rely on any spectral assumptions. In particular, it shows that the mixing time $t_\star$ is "optimal". It is based on the following lemma, which translates the fact that the geodesic flow has unit speed, and so the support of functions or measures does not grow too fast.
\begin{proposition}{}\label{prop:Croissancesupport}
	Consider $f \in L^2(M)$ such that $\Diam \Supp(f) \leq \delta$ for some $\delta > 0$. Then \[\Vol \Supp A_tf \leq Ce^{(d-1) t} \sinh\left((d-1)\delta\right),\]where $C>0$ depends only on the dimension $d$ of the manifold.	
\end{proposition}

\begin{proof}
	The situation can be lifted in $\mathds H^d$. Let us fix $p\colon \mathds H^d \mapsto M$ a covering map. Take a {fundamental Dirichlet domain (or any nice enough fundamental domain)} for the manifold $M$, and define a lift of $f_0$ by \[\forall x \in \mathds H^d, \quad \widetilde{f_0}(x) = \begin{cases} f_0(x)& \text{if $x \in \mathcal D$},\\0& \text{if $x \notin \mathcal D$.} \end{cases}\]

	Then $\widetilde{f_t} = A_t^{\mathds H^d}\widetilde{f}_0$ is a lift of $f_t$ in the sense that \[f_t(x) = \sum_{y\in p^{-1}(x)}^{} \widetilde{f_t}(y).\]
	
	In particular, since $p$ is a local isometry, $\Vol \Supp f_t \leq \ell_d \Supp \widetilde{f_t}$. When $f_0$ is supported inside a ball of radius $\delta$, then so is $\widetilde{f_0}$. As a consequence, $\widetilde{f_t}$ is supported in a $d$-dimensional spherical shell delimited by spheres of radiuses $t-\delta, t+\delta$. We recall that $\Gamma_d$ be the volume of the $d-1$-dimensional Euclidean sphere $\mathds S^{d-1}\subset \R^{d}$. Then the volume of any ball $B(y, r)$ of radius $r$ in $\mathds H^d$ is given by \[\Vol B(y, r) = \Gamma_d \int_0^{r}\sinh^{d-1} u\ \d u.\] As a consequence, \begin{align*} \ell_d \Supp \widetilde{f}_t \leq \Gamma_d \int_{t-\delta}^{t+\delta}\sinh^{d-1} u\ \d u &\leq {\Gamma_d}{} \int_{t-\delta}^{t+\delta}e^{(d-1)u} \ \d u\\ &\leq \frac{2\Gamma_d}{(d-1)} e^{(d-1) t} \sinh\left((d-1)\delta\right),
	\end{align*} and so  \[\Vol \Supp f_t \leq \frac{2\Gamma_d}{(d-1)} e^{(d-1) t} \sinh\left((d-1)\delta\right),\] as desired.
\end{proof}

This inequality is sufficient to get what we want:

\begin{proof}[Proof of Proposition \ref{prop:LB}]
	Consider the support $\Supp f_t$ of the density $f_t$. While there exists $0<\eta<1$ such that $\Vol \Supp f_t \leq \eta V_n$, we have  \[\DTV\left( \mu_{t}, \pi \right) \geq  \mu_t(\Supp f_t)-\pi(\Supp f_t) \geq 1 -\eta.\]
	According to Proposition $\ref{prop:Croissancesupport}$, \[\Vol \Supp f_t \leq Ce^{(d-1)t}\sinh\left((d-1) \delta\right).\] Fix $\varepsilon>0$ and take $t = (1-\varepsilon)t_{\star}$. Then
	\begin{equation*}\label{eq:suppft}
		\Vol \Supp f_{\left(1-\varepsilon\right)t_{\star}} \ll V_n,
	\end{equation*} where the last inequality holds in the three situations presented in Theorems \ref{th:D1}-\ref{th:S2}. It follows that \[\lim \DTV\left( \mu_{(1-\varepsilon)t_\star}, \pi \right) = 1,\] when $\delta \to 0$ for Theorem \ref{th:D1} and when $n\to\infty$ for Theorems \ref{th:S1} and \ref{th:S2}.
\end{proof}

As we will confirm in the next subsection, this lower bound is the correct estimate for the mixing time. This is, in fact, not very surprising. Indeed, when looking at the behavior of $\widetilde{f_{t_{\star}}}$ on $\mathds H^d$, we can already see that it is approximately uniform on its support. As a consequence, it only suffices that this support covers "quite uniformly" the surface to get the correct estimate. This will not be exactly the approach we choose to get the cutoff, but it is still an interesting heuristic.

\subsection{Upper bound} \label{subsec:UB}In this section, we conclude the proofs of Theorems \ref{th:D1}-\ref{th:S2} by showing the upper bound on the mixing time : \begin{itemize}
	\item Under the hypothesis of Theorem \ref{th:D1}, for any $\varepsilon >0$,
	\[\lim_{\delta \to 0}  \sup_{t \geq (1 + \eps) t_\star^\delta} \DTV (\mu^\delta_t,\pi) = 0.\]
	\item Under the hypothesis of Theorem \ref{th:S1} or Theorem \ref{th:S2}, for any $\varepsilon>0$,
	\[\lim_{n \to \infty}  \sup_{t \geq (1+\eps) t_\star^n} \DTV (\mu^n_t,\pi^n) = 0.\]
\end{itemize}

\subsubsection{Spectral decomposition of \texorpdfstring{$A_t$}{A\_t}}
The main step is to study $A_t$ from a spectral point of view. Unless otherwise stated, we will work on $L^2(M, \pi)$, abbreviated $L^2$, for the rest of the section. Using the fact that $\gamma_t(v)$ = $\gamma_{-t}(-v)$, it is not difficult to show that for any $t \geq 0$, $A_t$ is a self-adjoint operator which ensures that the spectrum is contained in $\R$. In what follows, we denote by $\boldsymbol{\Im(z)}$ the imaginary part of a complex number $z$. The aim of this part is to get the following propositions, which describe the spectrum of $A_t$ much more precisely. They will be at the core of the end of the proofs of Theorems \ref{th:D1}-\ref{th:S2}. In what follows, we denote by $\Im(z)$ the imaginary part of a complex number $z$.

\begin{proposition}{}\label{prop:SpAt0}
	Let $\phi_k$ be an eigenfunction of the Laplace-Beltrami operator, with eigenvalue $\lambda_k$. Then, $\phi_k$ is also an eigenfunction of $A_t$, and there exists a function $\nu_t\colon \Rp \to \R$ such that \[A_t\phi_k = \nu_t(\lambda_k)\phi_k.\]
	If $\lambda = \frac{(d-1)^2}{4}+u^2$, $u \in \R \cup i \left[-\frac{d-1}{2}, \frac{d-1}{2}\right]$,
	we have 
$$
\nu_t(\lambda) =  c_d \int_{-1}^1 ( \cosh t  +  x \sinh t  )^{iu - \frac{d-1}{2}}  ( 1 - x^2)^{ \frac{d - 3}{2}} \d x,  
$$
where $c_d = \frac{\Gamma\left( \frac{d}{2}\right)}{\sqrt{\pi} \Gamma\left( \frac{d-1}{2}\right)}$. In particular, for $d=3$ 
$$
\nu_t(\lambda) = \frac{\sin (ut)}{u \sinh (t)}.
$$
\end{proposition}

\begin{proposition}{}\label{prop:SpAt}
	Let $\lambda = \frac{(d-1)^2}{4}+u^2$ with $u \in \R_+ \cup i \left[0, \frac{d-1}{2}\right]$. Let $t_0 > 0$. There exists a constant $C = C(d,t_0)$ such that for all $t \geq t_0$, if $u \in i \left[0, \frac{d-1}{2}\right]$, we have
	$$|\nu_t(\lambda)|  \leq C  \min \left( t, \frac{1}{\Im(u)} \right)  e^{ \Im (u)  t-\frac{d-1}{2}t},$$
	and  if $u \in \R_+$,
	$$
	|\nu_t(\lambda)|  \leq C  \min \left( t, \frac{1}{u^{\frac{d-1}{2}}} \right)  e^{-\frac{d-1}{2}t}.
	$$
\end{proposition}

Note that we have given the closed-form expression for $\nu_t(\lambda)$ for $d=3$, but more generally for any odd $d \geq 3$, the integral in the definition of $\nu_t(\lambda)$ could be computed explicitly by integration by parts. For even $d \geq 2$, however, we do not know if it is possible to write a closed-form expression for $\nu_t(\lambda)$ (but see Proposition \ref{prop:SpAt} below for quantitative estimates).
Note also that in the case $d=2$, Proposition \ref{prop:SpAt} improves on the estimate made in \cite[Proposition 2.3]{Golubev} when $\lambda$ is large. This improvement will be crucial in the proof of Theorem \ref{th:D1}. 

Before proving the theorems, let us recall the spectrum of the Laplace-Beltrami operator. The expression of the Laplace-Beltrami operator $\Delta$ in coordinates for $C^\infty(M)$ functions is as follows: \[\forall f \in C^\infty(M), \quad \Delta f = - \frac{1}{\sqrt{g}}\sum_{i, j}^{}\partial_i\left(\sqrt{g}g^{ij}\partial_j f\right).\]

As an $L^2$-unbounded operator, $\Delta$ is symmetric, non-negative, and admits a unique self-adjoint extension (\cite{Stri}). The spectrum $\mathfrak S(\Delta)$ of the Laplacian operator (or of its extension) is well-understood. In the compact case, it is only pure point, and we have the following theorem, as reminded in the introduction: 
\begin{theorem}{(\cite{Chavel})}\label{th:LaplacianSpec}
Let $M$ be a compact connected Riemannian manifold without boundary. Then there exists a complete orthonormal system $(\phi_i)_i\geq 0$ of $L^2(M)$ composed of $C^\infty$-eigenfunctions of the Laplace-Beltrami operator with eigenvalues $0 = \lambda_0 < \lambda_1 \leq \lambda_2\leq \dots, \lambda_k \to \infty$ as $k\to \infty$.
\end{theorem}

	As a consequence, any $f \in L^2$ can be decomposed as \begin{equation}f = \sum_{\lambda_k \in \mathfrak S(\Delta)} \langle f, \phi_k \rangle \phi_k. \label{Dec}
	\end{equation} Combined with Proposition \ref{prop:SpAt}, this last expression is quite convenient for our problem, as we will see later on.
		
	We recall the notation $\sigma_d = \left(\frac{d-1}{2}\right)^2$ for the bottom of the $L^2$ spectrum of the Laplacian on $\dH^d$. As explained in \cite{Golubev} for the  special case of surfaces, the spectrum of the Laplace-Beltrami operator on the manifold $M$ can be parametrized by its unitary dual. We will not go into the details, but note that this leads to the following partition of the spectrum into three parts: \begin{itemize}
		\item the principal series containing all the eigenvalues $ \lambda_k = \left( \frac{d-1}{2}\right)^2 + u^2$ for some $u \in \R$, that is all the eigenvalues $\lambda_k \geq \sigma_d$. In the following, we write $\mathfrak S_P(\Delta)$ the set of eigenvalues of the Laplace-Beltrami operator coming from the principal series.
		\item the complementary series containing all the eigenvalues $\lambda_k = \left(\frac{d-1}{2}\right)^2 + u^2$ with $iu \in \left(-\frac{d-1}{2}, \frac{d-1}{2} \right) \setminus \left\{0\right\}$, that is $0 < \lambda_k < \sigma_d $. We write $\mathfrak S_C(\Delta)$ the eigenvalues coming from the complementary series.
		\item the discrete series corresponding to the eigenvalue $0 = \left(\frac{d-1}{2}\right)^2 + u^2$ for $iu = \pm \frac{d-1}{2}$, associated with the constant eigenfunctions. 
	\end{itemize} As we can see in Proposition \ref{prop:SpAt}, this trichotomy is of interest to us, because the behavior of the eigenvalues of $A_t$ depends strongly on which part of the spectrum they belong to. In particular, the behavior of the geodesic flow is controlled by how many eigenvalues belongs to the complementary series. Our assumptions about the spectrum exactly tell us that we can reasonably control their number. Note that choosing $u$ or $-u$ in the previous case does not change the associated eigenvalue.

Let us now prove Proposition \ref{prop:SpAt0}.	

\begin{proof}[Proof of Proposition \ref{prop:SpAt0}]
	It follows from the fact that $A_t$ is a radial operator.
	To prove it in the special case of $A_t$, we refer to \cite[Proposition A.1]{Boulanger}. In his case, Boulanger considers the semi-group of operators $\left(\mathcal{O}^{\rho}\right)_\rho$, where $\mathcal{O}^{\rho}$ does the mean on the hyperbolic ball of radius $\rho$ instead of the mean on the hyperbolic sphere of radius $\rho$. He proves what is known as Delsarte's formula, which states that \begin{equation}\label{Boul}\mathcal{O}^{\rho}(\phi_k) = \frac{1}{\Vol {B_o(\rho)} }\int_{B_o(\rho)}{ y^s}{ \frac{\d x_1\dots \d x_{d-1} \d y}{y^d}} \ \phi_k,\end{equation} where $s(d-1-s)  = \lambda_k = \lambda$, $o = (0, 0,\dots, 1)$ in the half-space model of $\mathds H^d$ and $B_o(t)$ is hyperbolic ball of radius $t$. Note that the value of the integral depends only on the original eigenvalue, and not on the eigenfunction nor the choice of $s$ between the two solutions of the equation.		
	Now, we write $\lambda_k = \frac{(d-1)^2}{4}+u_k^2$ and $s =s_k= \frac{d-1}{2} \pm i u_k$.
	
	A direct consequence of $\eqref{Boul}$ for us is that  \[A_t\phi_k =  \left(\frac{1}{\Vol {S_o(t)} }\left.\frac{\partial}{\partial \rho}\right\vert_{\rho = t}\int_{B_o(\rho)}{ y^s}{ \frac{\d x_1\dots \d x_{d-1} \d y}{y^d}} \right)\phi_k,\]
	where $S_o(t)$ is the $(d-1)$-dimensional hyperbolic sphere of radius $t$. We deduce that \[\nu_t(\lambda) = \frac{1}{\Vol {S_o(t)} }\left.\frac{\partial}{\partial \rho}\right\vert_{\rho = t}\int_{B_o(\rho)}{ y^s}{ \frac{\d x_1\dots \d x_{d-1} \d y}{y^d}}.\] 
	In the model of the hyperbolic half-space model, the ball $B_o(t)$ corresponds to the Euclidean ball of center $\cosh(t) \cdot o$ and radius $\sinh t$.  We use spherical coordinates centered in $\cosh t \cdot o$ \[\begin{cases}
		r &= d_h(x, o)\\
		x_1 &= \sinh r \cos\theta_1\\
		x_2 &= \sinh r \sin\theta_1\cos{\theta_2}\\
		\dots\\
		x_{d-1} &= \sinh r \sin\theta_1 \dots \cos \theta_{d-1}\\ 
		y &= \cosh r + \sinh r \sin\theta_1 \dots \sin \theta_{d-1}
	\end{cases} \quad \text{where } \begin{cases} r \in [0, t]\\ \theta_j \in [0, \pi]& \text{ for } j \leq d-2,\\ \theta_{d-1} \in (-\pi, \pi)\end{cases}\]
	to compute the integral. The Jacobian of this diffeomorphism is $y\sinh^{d-1} r \prod_{j}\sin^{d-1-j}\theta_j$. Let us write $\d\theta = \prod_{j= 1}^{d-1} \d\theta_j$ for the Lebesgue measure on $\R^{d-1}$. We find that  \begin{align*} \nu_t(\lambda) &=\frac{\sinh^{d-1} t}{\Vol S_o(t)} \int_{[0, \pi]^{d-2}\times(-\pi, \pi)}\left(\cosh t + \sinh t\prod_{j = 1}^{d-1} \sin \theta_j\right)^{s-d+1} \prod_{j = 1}^{d-1} \sin^{d-1-j}\theta_j \d\theta \\&=\frac{\sinh^{d-1} t}{\Vol S_o(t)} \int_{[-\frac{\pi}{2}, \frac{\pi}{2}]^{d-2}\times(-\pi, \pi)}\left(\cosh t + \sinh t \prod_{j=1}^{d-1} \cos \theta_j\right)^{iu-\frac{d-1}{2}} \prod_{j = 1}^{d-1} \cos^{d-1-j}\theta_j \d\theta.
	\end{align*} Because $\Vol S_o(t) = \Gamma_d \sinh^{d-1}t$ where $\Gamma_d =  \frac{2\pi^{d/2}}{ \Gamma (d/2)}$ is the volume of the $(d-1)$-dimensional Euclidean sphere, we get
	\begin{align*}
		\nu_t (\lambda) &=\frac{1}{\Gamma_d} \int_{(-\frac{\pi}{2}, \frac{\pi}{2})^{d-2}\times (-\pi, \pi)}\left(\cosh t + \sinh t\prod_{j=1}^{d-1} \cos \theta_j\right)^{iu-\frac{d-1}{2}} \prod_{j = 1}^{d-1} \cos^{d-1-j}\theta_j \d\theta \\ 
		&=  \frac{2^{d-1}}{\Gamma_d} \int_{[0, \frac{\pi}{2})^{d-2}\times[0, \pi)}\left(\cosh t + \sinh t\prod_{j=1}^{d-1} \cos \theta_j\right)^{iu-\frac{d-1}{2}} \prod_{j = 1}^{d-1} \cos^{d-1-j}\theta_j \d\theta.
	\end{align*}
	It remains to simplify the above integral. We perform a first change of variables from $[0, \frac{\pi}{2})^{d-2}\times[0, \pi)$ to $[0,1)^{d-2} \times (-1,1)$: $$ J (\theta_1,\dots,\theta_{d-1}) = (\cos (\theta_1), \dots,\cos(\theta_{d-1})).$$  We get: 
	\begin{align*}
		\nu_t (\lambda)  & = \frac{2^{d-1}}{\Gamma_d}  \int_{[0, 1)^{d-2}\times(-1, 1)}\left(\cosh t + \sinh t\prod_{j=1}^{d-1} x_j\right)^{iu-\frac{d-1}{2}}  \frac{ \prod_{j = 1}^{d-1} x_j^{d-1-j} \d x}{\prod_{j=1}^{d-1} \sqrt{1 - x_j^2}} \\
		&= \frac{2^{d-1}}{\Gamma_d} \left( I_+ + I_- \right),
	\end{align*}
	where $\d x = \prod_{j=1}^{d-1} \d x_j$ and
	$$
	I_{\pm} = \int_{[0, 1)^{d-1}}\left(\cosh t \pm \sinh t\prod_{j=1}^{d-1} x_j\right)^{iu-\frac{d-1}{2}}  \frac{ \prod_{j = 1}^{d-1} x_j^{d-1-j} \d x}{\prod_{j=1}^{d-1} \sqrt{1 - x_j^2}}.
	$$
	If $d =2$, it concludes the proof. We now assume $d\geq 3$. We perform the change of variable from $(0,1)^{d-1}$ to $\mathcal D = \left\{ x \in (0,1)^{d-1} : x_{d-1} < x_{d-2} \dots < x_1 \right\}$  defined by:
	$$
	(y_1,\dots,y_{d-1}) = K(x_1,\dots,x_{d-1}) = (x_1,x_1 x_2, \dots, x_{1} x_2 \dots x_{d-1}).
	$$
	Its inverse $K^{-1}$ is given by $K^{-1} (y)_j = y_{j}/y_{j-1}$ with the convention that $y_0 = 1$. In particular, $ D K^{-1} (y)$ is an upper triangular matrix and $\det ( D K^{-1} (y) ) = \prod_{j=1}^{d-2} y_{j}^{-1}$. Since $\prod_{j=1}^{d-1} x_j^{d-1-j}   = \prod_{j=1}^{d-2} y_j$, we deduce that, with $\d y = \prod_{j=1}^{d-1} \d y_j$,
	\begin{align*}
		I_{\pm} & = \int_{\mathcal D}\left(\cosh t \pm y_{d-1} \sinh t \right)^{iu-\frac{d-1}{2}}  \frac{   \d y}{\prod_{j=1}^{d-1} \sqrt{1 - y_j^2/y_{j-1}^2 }} \\
		& = \int_{\mathcal D}\left(\cosh t \pm y_{d-1} \sinh t \right)^{iu-\frac{d-1}{2}}  \frac{  \prod_{j=1}^{d-2} y_{j}  \d y}{\prod_{j=1}^{d-1} \sqrt{y_{j-1}^2  - y_j^2}} \\
		& =  2^{-d+2}  \int_{0}^1 \left(\cosh t \pm x \sinh t \right)^{iu-\frac{d-1}{2}} \varphi_d (x^2) \d x,
	\end{align*}
	with, for $x \geq 0$,
	$$
	\varphi_d (x) = \int_{\mathcal D_{d-2} (x)}  \frac{ \d z_1 \dots \d z_{d-2} }{\prod_{j=1}^{d-1} \sqrt{z_{j-1}  - z_j}} ,
	$$
	where we set $z_0 = 1$, $z_{d-1} = x$ and $\mathcal D_{d-2} (x) = \left\{ (z_1,\dots,z_{d-2}) : 1 > z_1 > \dots > z_{d-2} > x \right\}$ (we have performed the change of variable $z_j = y_j^2$ at the last step). 

	We check by recursion that $\varphi_{d} (x) = u_d (1-x)^{\frac{d-3}{2}}$ for some real $u_d >0$. For $d=2$, we have seen that $u_2= 1$. We use the recursion: 
	$$
	\varphi_{d+1} (x) = \int_{x}^1 \frac{\varphi_{d} (y)}{\sqrt{y - x}} \d y. 
	$$
	By the recursion hypothesis, we get
	\begin{align*}
		\varphi_{d+1} (1-x) & = \int_{0}^{x} \frac{\varphi_{d} (1-y)}{\sqrt{x-y}} \d y \\
		&  = u_d  \int_{0}^{x} \frac{y^{\frac{d-3}{2}}}{\sqrt{x-y}} \d y \\
		& = u_d  x^{\frac{d-2}{2}} \int_{0}^{1} \frac{(1-t)^{\frac{d-3}{2}}}{\sqrt{t}} \d t.  
	\end{align*}
	Hence $u_{d+1} = u_{d} \int_0^{1} (1-t)^{\frac{d-3}2} t^{-1/2} \d t $. We finally remark that 
	$$
	\int_{0}^{1} \frac{(1-t)^{k/2}}{\sqrt{t}} \d t = 2 \int_{0}^{1} (1-t^2)^{k/2} \d t = 2 \int_{0}^{\pi/2} \sin^{k+1} (\theta) \d\theta .
	$$
	We recognize the Wallis integral $W_{k+1}$ and the conclusion easily follows from the formulas \[W_k = \frac{\sqrt{\pi}}{2} \frac{\Gamma\left(\frac{k+1}{2}\right)}{\Gamma \left(\frac{k+2}{2}\right)} \quad \text{and} \quad \Gamma_d = \frac{2 \pi^{d/2}}{ \Gamma\left(d/2\right)}.\] The integral formula for $\nu_t(\lambda)$ follows. For $d=3$, it is obvious to integrate. \end{proof}

Proposition \ref{prop:SpAt0} can be used to obtain the quantitative estimates on $\nu_t(\lambda)$ of Proposition \ref{prop:SpAt}.

\begin{proof}[Proof of Proposition \ref{prop:SpAt}]
	For $d=3$, the lemma is clear from the explicit expression given in Proposition \ref{prop:SpAt0}. For general $d \geq 2$, we start by proving that all $u \in \R_+ \cup i \left[0, \frac{d-1}{2}\right]$,
	\begin{equation}\label{eq:bdnut1}
		|\nu_t(\lambda)|  \leq C  \min \left( t, \frac{1}{\Im(u)} \right)  e^{ \Im (u)  t-\frac{d-1}{2}t}.
	\end{equation}
	Let $v = \Im(u) \geq 0$ and $a = \tanh t \in [\tanh t_0,1)$. By Proposition \ref{prop:SpAt0} (applied to $-u$), 
	\begin{align*}
		|  \nu_t (\lambda) | &= c_d (\cosh t )^{v - \frac{d-1}{2}} \left| \int_{-1}^{1} ( 1+ a x)^{-iu - \frac{d-1}{2}} ( 1- x^2)^{ \frac{d - 3}{2}} \d x \right| \\
		& \leq c_d (\cosh t )^{v - \frac{d-1}{2}}  \int_{-1}^{1} ( 1+ a x)^{v - \frac{d-1}{2}} ( 1- x^2)^{ \frac{d - 3}{2}}\d x \\
		& \leq c_d e^{v  t-\frac{d-1}{2}t} \left( I_+ + I_-\right),
	\end{align*}
	where 
	$$
	I_{\pm} = \int_{0}^{1} ( 1 \pm a x)^{v - \frac{d-1}{2}} ( 1- x^2)^{ \frac{d - 3}{2}}\d x.
	$$
	We have $I_+ \leq  \int_{0}^{1}   ( 1- x^2)^{ \frac{d - 3}{2}}\d x < \infty$. On the other hand, setting $a =  1- \delta$ and $b = a / \delta = (1 -\delta) / \delta$, we get through the change of variable $x \to 1 -x$:
	\begin{align*}
		I_- & = 	\int_{0}^{1} ( 1 - a (1-x) )^{v - \frac{d-1}{2}} ( x (2 -x) )^{ \frac{d - 3}{2}}\d x	\\
		& =  \delta^{v - \frac{d-1}{2}} \int_{0}^{1} ( 1 + b x )^{v - \frac{d-1}{2}} ( x (2 -x) )^{ \frac{d - 3}{2}}\d x \\
		& \leq C \delta^{v - \frac{d-1}{2}} \int_{0}^{1} ( 1 + b x )^{v - \frac{d-1}{2}}  x  ^{ \frac{d - 3}{2}}\d x, 
	\end{align*}		
	where at the last line, we have used that for all $x \in [0,1]$, $(2 -x)^{ \frac{d - 3}{2}} \leq C$ with $C = 1$ for $d \in \{2,3\}$ and $C = 2^{(d-3)/2}$ for $d \geq 4$. With the change of variable $y = bx$ and recalling that $\delta b = a \in [\tanh t_0,1)$, we find: 
	$$
	I_- \leq C a^{-\frac{d-1}{2}} \delta^{v}  \int_{0}^{b} ( 1 + y )^{v - \frac{d-1}{2}}  y  ^{ \frac{d - 3}{2}}\d y.
	$$
	We have $b = b(t) = (e^{2t} - 1)/ 2 \in [b_0,e^{2t}]$ with $b_0 = b(t_0) >0$. We write
	\begin{align*}
		\delta^{v}  \int_{0}^{b} ( 1 + y )^{v - \frac{d-1}{2}}  y  ^{ \frac{d - 3}{2}}\d y & = \delta^{v}  \int_{0}^{b_0} \dots +  \delta^{v}  \int_{b_0}^{b} \dots.
	\end{align*}
	Note that in a neighborhood of $0$, the integrand is equivalent to $y  ^{ \frac{d - 3}{2}}$, hence 
	$$
	\delta^{v}  \int_{0}^{b_0} ( 1 + y )^{v - \frac{d-1}{2}}  y  ^{ \frac{d - 3}{2}}\d y \leq C
	$$
	for some constant $C = C(d,t_0)$ (recall that $v \leq (d-1)/2$). Similarly, 
	$$
	\delta^{v}  \int_{b_0}^{b} ( 1 + y )^{v - \frac{d-1}{2}}  y  ^{ \frac{d - 3}{2}}\d y \leq \delta^{v}  \int_{b_0}^{b} ( 1 + y )^{v} y  ^{ -1} \d y \leq C \delta^{v}  \int_{b_0}^{b} y  ^{v -1} \d y,
	$$
	for some constant $C = C(d,t_0)$.  We finally use that
	$$
	\int_{b_0}^{b} y  ^{v -1} \d y  = \frac{1}{v}(b^v - b_0^v)= \frac{b^v}{v}\left(1 - \frac{b_0^v}{b^v}\right) \leq b^v \min \left( \frac{1}{v} , \ln \frac{b}{b_0}\right),
	$$
	Using again $\delta b = a \leq 1$ and $\ln (b) \sim 2t$ as $t \to \infty$, this proves \eqref{eq:bdnut1}.
	
	It remains to prove that if $u \in \R_+$, we have
	\begin{equation}\label{eq:bdnut2}
		|\nu_t(\lambda)|  \leq   \frac{C}{u^{\frac{d-1}{2}}}   e^{-\frac{d-1}{2}t}. 
	\end{equation}
	This bound improves  \eqref{eq:bdnut1} for large $u$. It is obtained by taking into account the oscillatory part in the formula for $\nu_t(\lambda)$. With $ a = 1- \delta = \tanh (t)$ as above, from Proposition \ref{prop:SpAt0} we write
	\begin{align*}
		\nu_t (\lambda)  &= c_d (\cosh t )^{iu - \frac{d-1}{2}}   \int_{-1}^{1} ( 1+ a x)^{iu - \frac{d-1}{2}} ( 1- x^2)^{ \frac{d - 3}{2}} \d x.
	\end{align*}

For $d \geq 4$, we start by performing integration by parts to take advantage of the factor $( 1- x^2)^{ \frac{d - 3}{2}}$. 
To this end, we fix $\beta \in \mathbb{C}$ and $\alpha \in [0,1)$. For integers $ n,k \geq 0$ with $2k \leq n$, let $\mathcal P_{k,n} \subset \mathbb{C}[X]$ be the set of polynomials of degree at most $n$ which are divisible by $(1-x^2)^k$.  For $P \in \mathbb{C}[X]$, we define
$$
I_k (P) = \int_{-1}^{1} ( 1+ a x)^{\beta - k-1} P(x) ( 1- x^2)^{ -\alpha} dx.
$$
By integration by parts, if $k \geq 1$ and $P \in \mathcal P_{k,n}$, we observe that 
\begin{align*}
I_k (P) & = \frac{1} {a (\beta -k)}   \int_{-1}^{1} ( 1+ a x)^{\beta - k}  \left (  -2 \alpha \frac{ x P(x)}{1-x^2} - P'(x) \right) ( 1- x^2)^{ -\alpha} dx \\
& = \frac{1} {a (\beta -k)} I_{k-1} (Q_1),
\end{align*}
	  where $Q_1 \in \mathcal P_{k-1,n-1}$. In particular, by iteration, we deduce that there exists a polynomial $Q_k\in \mathcal P_{0,n-k}$ such that 
$$
I_k (P) =  \frac{ I_0(Q_k)} { a^k  \prod_{j=0}^{k-1}  (\beta - k  + j) } . 
$$
Note that $Q_k$ does not depend on the parameters $(a,\beta)$.
	  We apply this observation to $k = \lceil (d-3)/2 \rceil$, $P(x) = (1-x^2)^{   k  }$, $\alpha = 0$ for odd $d \geq 5$, $\alpha = 1/2$ for even $d \geq 4$ and  $\beta = iu + \alpha$ (so that $\frac{d-3}{2} = k - \alpha$, $\frac{d-1}{2} = k - \alpha+1 $). We deduce that for some polynomial $Q$ of degree  at most $k$ we have:
	  				\begin{align*}
	  \nu_t (\lambda)  &= \frac{ c_d (\cosh t )^{iu - \frac{d-1}{2}}}  { a^k  \prod_{j=0}^{k-1}  (i u - \frac{d-3}{2} +j ) }  \int_{-1}^{1} ( 1+ a x)^{iu - 1 + \alpha} Q(x) (1-x^2)^{-\alpha}dx \\
	  & =  \frac{ c_d (\cosh t )^{iu - \frac{d-1}{2}}}  { a^k  \prod_{j=0}^{k-1}  (i u -  \frac{d-3}{2} +j ) } \left( J_+ + J_-\right),  
	   \end{align*}
		with
		$$
		J_{\pm} = \int_{0}^{1} ( 1 \pm a x)^{iu - 1 + \alpha} Q(x) ( 1- x^2)^{ -\alpha}dx.
		$$	
		The proof of \eqref{eq:bdnut2} will be complete if we prove that 
\begin{equation}\label{eq:bdnut2b}
		|J_{\pm}|  \leq   \frac{C}{u^{1-\alpha}}. 
\end{equation}
	We start with $J_-$ and set $b = a /  \delta$ as above. Since $((1-x)^l)_{l \geq 0}$ is a basis of $\C[X]$, by linearity, we may assume that $Q(x) = (1-x)^l$ for some integer $0 \leq l \leq k$. From what precedes, we get
	\begin{align*}
		J_- &=  \delta^{iu - 1 + \alpha} \int_{0}^{1} ( 1+ b x )^{iu - 1 + \alpha} x^l ( x(2- x))^{ -\alpha } \d x  \\
		& = \delta^{iu} u^{-1} (\delta b)^{-1 + \alpha } b ^{-l} \int_{0}^{u \ln (1 + b)} e^{iy } e^{y \alpha / u }  (e^{y/u} -1)^{l - \alpha}  (2-(e^{y/u} -1)/b )^{- \alpha } \d y,
	\end{align*}
	where we have performed the change of variable $y = u \ln ( 1 + bx)$, that is $x = (e^{y/u} -1)/b$. Since $\delta b = a$, we find
	\begin{align*}
		J_- 
		& = \delta^{iu} u^{-1+\alpha} a^{-1 + \alpha } b ^{-l} \int_{0}^{u \ln (1 + b)} e^{iy} y^{-\alpha} f_b\left(\frac y u \right) \d y,
	\end{align*}
	where we have set 
	$$
	f_b(x) =  \frac{(e^{x} -1)^{l} }{ \left( g(x)  h_b(x) \right)^{\alpha} }, \quad g(x) =  \frac{1-e^{-x} }{x} , \quad h_b(x) =  2 - \frac{e^{x} -1}{b}.
	$$		
	The functions $g,h_b$ and $x \mapsto (e^x-1)^l$ are analytic on $\C$. If $\alpha = 0$, we deduce that $f_b(x)$ is analytic on $\C$. If $\alpha = 1/2$, consider the principal branch of the function $z \mapsto z^{-\alpha}$ analytic on $\C \backslash \R^-$ with $\R^- = (-\infty,0]$.  For $\theta \geq 0$, let $C_\theta = \{ x \in \C : \mathrm{arg} (x) \in [-\theta,\theta] \}$ where $\mathrm{arg}(x) \in [-\pi,\pi)$ is the argument of a complex number $x$. Let $C^*_\theta = \C_\theta \backslash \{0\}$ and $C_\theta (r) = C_\theta \cap B(r)$ where $B(r) = \{ x : |x| \leq r\}$. We claim that the image of $C_{\pi/8}  (\ln (1+b))$ by $x \mapsto g(x)  h_b(x)$ is contained in $C^*_{7\pi / 8}$. This claim implies that  $f_b$ is analytic on $C_{\pi/8}  (\ln (1+b))$.

Let us prove the claim. First, for all $x \in C_{\pi/2} ( \ln(1+b))$, it is easy to check that  $|e^x - 1|/b \leq (e^{|x|} - 1 ) / b \leq 1$. Thus, for all $x \in C_{\pi/2} ( \ln(1+b))$, $h_b(x)$ is in a ball of radius $1$ centered at $2$. In particular, we have $h_b(x) \in C^*_{\pi/4}$. Note also that $g(0) = 1$ and $g(x) \ne 0$ for all $x \in \C$. Also, for $ x \in C^*_\theta$ with $\theta \leq \pi/2$, we have $|\mathrm{arg} (g(x))| \leq |\mathrm{arg} (1 - e^{-x})| + |\mathrm{arg} (x^{-1})| \leq \pi / 2 + \theta$ because $\Re(1 - e^{-x}) = 1 -e^{-|x| \cos (\theta)} \cos (|x| \sin (\theta))  >0$. So we have checked that for all $x \in C^*_\theta$, we have $g(x)  h_b(x) \in C^*_{\pi/2 + \pi/4 + \theta}$. Since $g(0)  h_b(0) = 2$, the latter is true for all $x \in C_\theta$. We deduce the claim.

We can now complete the proof of \eqref{eq:bdnut2b} for $J_-$. Since  $f_b$ is analytic on $C_{\pi/8}  (\ln (1+b))$, we can apply Cauchy formula around a counterclockwise contour around the set $[0,u \ln (1+b)] \cup \{ x : x =  u \ln (1+b) e^{i \theta} , \theta \in [0,\pi/8]\} \cup e^{i \pi/8} [0,u \ln (1+b)]$. We get, with $r = u \ln (1+ b)$,
\begin{align*}
	K_- & = \int_{0}^{r} e^{iy} y^{-\alpha} f_b\left(\frac y u \right) \d y \\
	& = \int_{0}^{r} e^{iy  e^{i\pi/8}} ( e^{i\pi/8} y)^{-\alpha} f_b\left(\frac {y  e^{i\pi/8}}{ u}  \right) e^{i\pi/8}\d y -  \int_{0}^{\pi/8} e^{iy e^{i \theta} } (r e^{i\theta}) ^{-\alpha} f_b\left(\frac {r e^{i \theta}}{ u}  \right) r \d \theta.
\end{align*}
Taking absolute values, with $c = \sin (\pi/8) > 0$, we get 
$$
| K_- | \leq \int_{0}^{r} e^{-c y } y^{-\alpha} \left| f_b\left(\frac {y  e^{i\pi/8}}{ u}  \right)\right|  \d y + r^{1-\alpha} \int_{0}^{\pi/8} e^{- r \sin ( \theta) }  \left|f_b\left(\frac {r e^{i \theta}}{ u}  \right)\right|  \d \theta.
$$
Next we can use $\sin(\theta ) \geq 2 \theta / \pi$ and that for all $x \in C_{\pi/2}( \ln(1+b))$ we have $|f_b(x)| \leq C b^l$  for some constant $C > 0$. We get that
$$
| K_- | \leq C  b^l \int_0^\infty e^{-cy} y^{-\alpha} \d y + C b^l r^{1-\alpha}  \int_{0}^{\pi/8} e^{- 2 r \theta / \pi } \d\theta \leq C' b^l ( 1+ r^{-\alpha}),  
$$
for some new constant $C'$. Since $J_- = \delta^{iu} u^{-1+\alpha} a^{-1 + \alpha } b ^{-l} K_-$, this completes the proof of \eqref{eq:bdnut2b} for $J_-$.

It remains to prove \eqref{eq:bdnut2} for $J_+$. The proof is the same as the proof for $J_-$ so we will use the same notation. Using the same argument as above, we assume that $Q(x) = (1-x)^l$. We now set  $\delta = (1+a)$ and $b = a / (1 + a) \in [b_0,1/2)$ with $b _0 = b(t_0) \in (0,1/2)$. We write 
\begin{align*}
	J_+ &=  \delta^{iu - 1 + \alpha} \int_{0}^{1} ( 1 - b x )^{iu - 1 + \alpha} x^l ( x(2- x))^{ -\alpha } \d x  \\
	& = \delta^{iu} u^{-1} (\delta b)^{-1 + \alpha } b ^{-l} \int_{0}^{- u \ln (1 - b)} e^{-iy } e^{-y \alpha / u }  (1-e^{-y/u})^{l - \alpha}  (2-(1 - e^{-y/u})/b )^{- \alpha } \d y,
\end{align*}
where we have performed the change of variable $y =  - u \ln ( 1 - bx)$, that is $x = (1 - e^{-y/u})/b$. Since $\delta b = a$, we find
\begin{align*}
	J_+ 
	& = \delta^{iu} u^{-1+\alpha} a^{-1 + \alpha } b ^{-l} \int_{0}^{- u \ln (1 - b)} e^{-iy} y^{-\alpha} f_b\left(\frac y u \right) \d y,
\end{align*}
where now we have
$$
f_b(x) =  \frac{(1-e^{-x})^{l} }{ \left( g(x)  h_b(x) \right)^{\alpha} }, \quad g(x) =  \frac{e^{x} -1 }{x} , \quad h_b(x) =  2 - \frac{1- e^{-x}}{b}.
$$		
We can then conclude similarly thanks to a contour integration around the set $[0, - u \ln (1 - b)] \cup \{ x : x =  -u \ln (1+b) e^{ - i \theta} , \theta \in [0,\pi/8]\} \cup e^{ - i \pi/8} [0, - u \ln (1-b)]$ (note that in this case, it is easier since $\ln (1-b) $ and $b$ are uniformly bounded by some $C = C(t_0)$). 
\end{proof}

\subsubsection{Spectral bound on total variation distance}
In this section, we give a spectral bound on the total variation distance using Proposition \ref{prop:SpAt} and Equation $\eqref{Dec}$. This, combined with Proposition $\ref{prop:LB}$, is essentially sufficient to conclude the proofs of Theorems $\ref{th:D1}-\ref{th:S2}$. Then, although the outline of the proof is the same for each of our three frameworks, there are some subtleties when going from one to the other. We will therefore separate the three proofs into three subsections, starting with the case of Theorem \ref{th:S1}, which is the most direct. We will then deal with Theorem \ref{th:D1}, for which we need a little refinement, and conclude with the proof of Theorem \ref{th:S2}, which is the more computational situation.

In the following, for $\lambda_k \in \mathfrak S_C(\Delta)$, we define $s_k\in (0,1)$ as the only positive number such that $$\lambda_k = \sigma_d (1 -s^2_k)$$ (\textit{i.e.} $s_k  = \frac{2 iu}{d-1}$, where $u$ has been chosen as in the proof of Proposition \ref{prop:SpAt}, $\Im u < 0$), so that the conclusion of Proposition \ref{prop:SpAt} can be conveniently rewritten as:
\begin{equation} \label{eq:nulkb}
	|\nu_t(\lambda_k)| = \begin{cases}  O\left(t e^{-\frac{d-1}{2} t}\right)&\text{if $\lambda_k  \geq \sigma_d \in \mathfrak S_P $}\\O\left( t e^{-\frac{d-1}{2} (1 - s_k) t }\right)&\text{if $\lambda_k = \sigma_d ( 1 -s_{k}^2) \in \mathfrak S_C$}.\end{cases}
\end{equation}

\begin{proposition}{}\label{prop:MajL2}
	There exists a numerical constant $C >0$ such that for any $f \in L^2(M, \pi)$, we have \[\lVert A_t f-1\rVert_{L^2}^2 \leq C t ^2 \left(\sum_{\lambda_k \in \mathfrak S_C(\Delta)} e^{-(1-s_k)(d-1)t}\langle f_0, \phi_k \rangle^2 + \sum_{\lambda_k \in \mathfrak S_P(\Delta)} e^{-(d-1)t}\langle f_0, \phi_k \rangle^2 \right).         \]
	
	In particular, \[\DTV\left( \mu_t,\pi \right)^2 \leq C t^2 \left(\sum_{\lambda_k \in \mathfrak S_C(\Delta)} e^{-(1-s_k)(d-1)t}\langle f_0, \phi_k \rangle^2 + \sum_{\lambda_k \in \mathfrak S_P(\Delta)} e^{-(d-1)t}\langle f_0, \phi_k \rangle^2\right) .\]
\end{proposition} 
\begin{proof}
	
	By the Cauchy-Schwarz inequality,
	\begin{align*}
		\DTV\left( \mu_t,\pi \right) = \frac{1}{2}\lVert f_t - 1\rVert_{L^1(\pi)} &\leq \frac{1}{2} \lVert f_t-1\rVert_{L^2(\pi)}
		\\&\leq \frac{1}{2} \lVert A_t \cdot \left(f_0-1\right)\rVert_{L^2(\pi)}.
	\end{align*}
	On the other hand, the decomposition $\eqref{Dec}$ leads to
	\begin{align*}\lVert A_t\left(f_0-1\right) \rVert_{L^2(\pi)}^2 &=  \sum_{\lambda_k \in \mathfrak S_C(\Delta)} \nu_t(\lambda_k)^2\langle f_0, \phi_k \rangle^2 + \sum_{\lambda_k \in \mathfrak S_P(\Delta)} \nu_t(\lambda_k)^2\langle f_0, \phi_k \rangle^2 .
	\end{align*}
	It remains to use the upper bound \eqref{eq:nulkb}.
\end{proof}

\subsubsection{Proof of Theorem  \ref{th:S1}}

Theorem \ref{th:S1} is a simple consequence of Proposition \ref{prop:MajL2}. Indeed, it tells that
\begin{align*}\DTV\left( \mu_t,\pi \right)^2 &\leq C t \left(\sum_{\lambda_k \in \mathfrak S_C(\Delta)} e^{-(1-s_k) (d-1) t}\langle f_0, \phi_k \rangle^2 + \sum_{\lambda_k \in \mathfrak S_P(\Delta)} e^{-(d-1)t}\langle f_0, \phi_k \rangle^2\right)
	\\&\leq C t  e^{ -(1-s_{1}) (d-1) t} \lVert f_0\rVert_2^2
	\\& \leq  C t \frac{V_n}{\Vol B ( x_n, \delta)}  e^{ - (1-s_{1}) (d-1) t}.
\end{align*}

Fix $\varepsilon > 0$, and take $t \geq \left(1+\varepsilon \right)t^n_\star = \left(1+\varepsilon \right) \frac{\ln V_n}{d-1}$. Then, \[\DTV\left( \mu_t,\pi \right)^2 \leq {\frac{C t }{\Vol B ( x_n , \delta)}}V_n^{s_{1}(1+\varepsilon)-\varepsilon} \leq  {\frac{C' }{\Vol B ( x_n , \delta)}}V_n^{s_{1}(1+\varepsilon)-\varepsilon/2} ,\]
where $C'$ depends on $\varepsilon$ is such that $ C t \leq C' \exp( t (d-1) \varepsilon / 2)$ for all $t\geq 0$. 

By assumption, when $n$ becomes large enough, $s_{1}$ becomes small, so $s_{1}(1+\varepsilon) < \varepsilon / 4$. Since $\sinh x \geq x$, we get that $\Vol B ( x_n , \delta ) \geq C_d \min(\delta^d, \rho(x_n)^d)$ for some constant $C_d$ which depends only on the dimension $d$.  Therefore, \[\DTV\left( \mu_t,\pi \right)^2 \leq {\frac{C' }{ C_d \min(\delta,\rho(x_n))^d}}{V_n^{-\varepsilon/4}}.\]

If $\liminf_n \frac{\ln  \rho(x_n)}{\ln V_n} = 0$, then \[\frac{\varepsilon}{4} \ln V_n + d\ln \rho = \frac{\varepsilon}{4} \ln V_n + d \frac{\ln \rho}{\ln V_n} \ln V_n \limitninf{}{} \infty,\] and so $\DTV\left( \mu_t,\pi \right)^2 \limitninf{}{} 0$ as required. 	  \qed 

Note that using the improved estimates given in \cite{HideMaceraThomas} for $\lambda_1(S_g)$ when $S_g$ is sampled according to the Weil-Petersson distribution, one can obtain from the previous computation a windows of order at most $\ln \ln n$ for the cutoff for typical hyperbolic surfaces.

\subsubsection{Proof of Theorem \ref{th:D1}}

The proof of Theorem \ref{th:D1} follows the same lines but is slightly more subtle: at time 0, the density is too concentrated, leading to an $L^2$-norm in $\delta^{-d}$ for $f_0$. However, we use the fact that the mass of $f_0$ is concentrated on high eigenfunctions and take advantage of our estimate on $\nu_t(\lambda)$ for large $\lambda$ in Proposition \ref{prop:SpAt}.

We fix $\varepsilon > 0$ and take 	$t \geq (1+\varepsilon)t_\star$ with  $t_\star  = t_\star^\delta = - \frac{\ln \delta}{d-1}$. Arguing as in the proof of Proposition \ref{prop:MajL2},
\begin{align*} \DTV\left( \mu_t,\pi \right)^2 
	&\leq \frac{1}{4} \lVert A_t(f_0 - 1)\rVert^2_{L^2}			
	\\&= \frac{1}{4}\lVert \sum_{k \geq 1} \nu_t(\lambda_k)\langle f_0, \phi_k \rangle \phi_k \rVert^2_{L^2} \\
	& = \frac 1 4 \sum_{k \geq 1} \nu_t(\lambda_k)^2 \langle f_0, \phi_k \rangle^2 ,
\end{align*}
where we have used the orthogonality of the basis $(\phi_k)$ at the last line. 	Fix $k_0 \geq 1$ large enough so that $\lambda_k \geq 2 \sigma_d$ for all $k \geq k_0$. Let $k_\delta > k_0$ be defined later. We decompose the sum into three parts, 
\begin{equation*}
	\sum_{k \geq 1} \nu_t(\lambda_k)^2 \langle f_0, \phi_k \rangle^2 = \sum_{1 \leq k \leq k_0} \dots + \sum_{ k_0 < k < k_\delta} \dots +   \sum_{ k \geq k_\delta} \dots = S_{1}(t) + S_{2} (t) + S_{3} (t).
\end{equation*}
Recall that $f_0$ (and hence $S_i$, $i \in \{1,2,3\}$) depends implicitly on $\delta$ and the initial condition $x \in M$. 
Our goal is to prove that for all initial conditions $x \in M$, we have for $i \in \{1,2,3\}$,  		
\begin{equation}\label{eq:limSi}
\lim_{\delta \to 0} \sup_{t \geq (1+\varepsilon) t_{\star}}  S_{i}(t) = 0.
\end{equation}

To this end, we will rely on two classical results in differential geometry. First, the Weyl law asserts that 
\begin{equation}\label{eq:weyl}
\lambda_k \sim c_1 k^{2/d} \; \hbox{ as $k \to \infty$ for some $c_1 = c_1(d) >0$}.
\end{equation}
Moreover, the H\"ormander Theorem implies that as $l \to \infty$,
$$
\sup_{y \in M} \sum_{k=0}^l \phi^2_k(y) \sim l, 
$$
see \cite[Theorem 16.1]{zbMATH01637290}. Hence, for all $l \geq 1$,  $\sup_{y \in M} \sum_{k=1}^l \phi_k(y)^2 \leq c_2 l$ for some $c_2 = c_2(M)$, and we find, from Cauchy-Schwarz inequality,  
\begin{equation}\label{eq:lwl}
 \sum_{k=1}^l \langle f_0, \phi_k \rangle^2  \leq \sum_{k=1}^l \int_M \phi^2_k(y) f_0(y) d\pi(y) \leq c_2 l.
\end{equation}

In \eqref{eq:limSi}, we start with $S_1$.  Since $\lambda_1 > 0$,  Proposition \ref{prop:SpAt} implies that for any $k \geq 1$, we have $\nu_t (\lambda_k)^2  \leq C \exp ( -  \alpha t)$ for some $C,\alpha >0$. We deduce that, with $C' = c_2 k_0 C$,
$$
S_{1}(t) \leq C' e^{- \alpha t}.
$$
The right hand does not depend on the initial condition $x \in M$ and $\delta > 0$. Since $t \to \infty$, we deduce that \eqref{eq:limSi} holds for $i =1$.

We now deal with $S_2(t)$ and $S_3(t)$. From Proposition \ref{prop:SpAt},  for all $\lambda \geq 2 \sigma_d$ and $t \geq t_0  >0 $, we have 
$$
| \nu_t (\lambda) | \leq   C  \lambda^{-\frac{d-1} 4} e^{- \frac{d-1}{2} t },
$$
where $C = C(t_0)$ (we have used that $\lambda = \sigma_d + u^2$ is asymptotically equivalent to $u^2$ for all large $\lambda$). Hence by Weyl law \eqref{eq:weyl}, for some new $C >0$, for all $k \geq k_0$, 
\begin{equation}\label{eq:weyl2}
| \nu_t (\lambda_k) |^2 \leq   C  k^{- 1 + 1/d} e^{- (d-1) t }.
\end{equation}

We deduce the upper bound: 
$$
S_{3}(t)  \leq  C e^{- (d-1) t}  k_\delta^{- 1 + 1/d} \sum_{k \geq k_\delta} \langle f_0 , \varphi_k \rangle^2 \leq C e^{- (d-1) t}  k_\delta^{- 1 + 1/d} \| f_0 \|^2_{L^2} .
$$
%We have $\| f_0 \|^2_{L^2}  = O(\delta^{-d})$. Also by Weyl's law, $\lambda_k \sim c k^{2/d}$ as $k \to \infty$ for some $c >0$. Hence, for some new $C >0$, for all $k \geq k_\delta$,
%$$
%| \nu_t (\lambda_k) |^2 \leq   C  k^{- 1 + 1/d} e^{- (d-1) t },
%$$
We have $\| f_0 \|^2_{L^2}  = O(\delta^{-d})$. Picking $k_\delta = \lceil \delta^{-d} \rceil \vee (k_0+1)$, we get, for some new $C >0$, 
$$
S_{3}(t)  \leq C \delta^{-1} e^{- (d-1) t}.
$$
We deduce that for any $T > 0$, for all $t \geq t_\star + T$, 
$S_{3}(t) \leq C e^{- T}$.   In particular, \eqref{eq:limSi} holds for $i =3$.

Similarly, using  \eqref{eq:weyl2}, 
$$
S_{2} (t) \leq C e^{- (d-1) t} \sum_{ 1 \leq k <  k_\delta} \frac{ \langle f_0 , \phi_k \rangle^2 }{k^{1 - 1/d}}. 
$$
For integer $l \geq 1$, we set $N_y (l) = \sum_{1 \leq k \leq l } \langle f_0 , \phi_k(y) \rangle^2$. By \eqref{eq:lwl}, $N_y (l) \leq c_2 l$. Hence, by discrete integration by parts,
\begin{align*}
\sum_{ 1 \leq k < k_\delta} \frac{\phi^2_k(y) }{k^{1 - 1/d}} & =  \frac{N_y(k_\delta-1)}{k_\delta^{1-1/d}} + \sum_{1 \leq k < k_\delta} N_y(k) \left(\frac{1}{k^{1-1/d}} - \frac{1}{(k+1)^{1-1/d}}\right) \\
& \leq  C' \left( k_\delta^{1/d} + \sum_{1 \leq k < k_\delta}  \frac{1}{k^{1-1/d}} \right) \\
& \leq C''k_\delta^{1/d}. 
\end{align*}
For our choice  $k_\delta = \lceil \delta^{-d} \rceil \vee (k_0+1)$, we deduce that, for some new $C >0$, 
$$
S_2 (t) \leq C \delta^{-1} e^{- (d-1) t}.
$$
Arguing as above,  \eqref{eq:limSi} holds for $i =2$. It concludes the proof.
\qed

\subsubsection{Proof of Theorem \ref{th:S2}}

We first remark that on $T_n$, we have that $\inf_{x \in T_n} \Vol B (x , \delta) \geq C_d \min(\delta, \rho(x))^d \geq  V_n^{-\eta_n}$ for some sequence $\eta_n \to 0$. We condition on $T_n$ and prove that there exists a sequence $S_n \subset T_n$ such that $\lim_n \pi_n(S_n\vert T_n) =1$ and for any sequence $x_n \in S_n$, the process $X^n_{x_n,t}$ exhibits cutoff at time $t_\star = t^n_\star = \frac{\ln V_n}{d-1}$ when $n \to \infty$. Because $\pi(T_n) \to 1$, this is sufficient. 			We denote $\widetilde{\pi} = \pi(. |T_n)$.  

Fix $\varepsilon >0$, and take $t \geq (1+\varepsilon)t_\star$. 			Let us estimate $\DTV\left( \mu_t, \pi\right)$ when the initial condition is centered at $x$ sampled according to $\widetilde \pi$. As before, we start from Proposition \ref{prop:MajL2} to get that \[\DTV\left( \mu_t,\pi \right)^2 \leq C t \left( \sum_{\lambda_k \in \mathfrak S_C} e^{-(1-s_k)(d-1)t}\langle f_0, \phi_k \rangle^2 + \sum_{\lambda_k \in \mathfrak S_P} e^{-(d-1)t}\langle f_0, \phi_k \rangle^2 \right).\]
Note that since $\liminf_n \lambda_1(M) > 0$, we have $\limsup_n s_1 < 1$ and the right-hand side of the above expression is decreasing in $t$ for $t$ large enough (since $t e^{-\alpha t}$ is decreasing for $t \geq 1/\alpha$).  			As a consequence, remembering that $\mu_t$ and $f_0$ depend on the point $x$ where the ball is centered, we get that, for $x$ sampled according to $\widetilde \pi$, 
\begin{eqnarray}\label{eq:Etilde}
\E_{\widetilde\pi} \sup_{s \geq t} \DTV\left( \mu_s, \pi\right)^2  \leq Ct e^{-(d-1)t} \left(\sum_{\lambda_k \in \mathfrak S_C} e^{  s_k (d-1)t}\int_{x \in M} \langle f_0, \phi_k \rangle^2  \d\widetilde\pi(x) +  \int_{x \in M} \lVert f_0 \rVert_2^2 \d\widetilde\pi(x)\right).
\end{eqnarray}

For $n$ large enough, $\pi(T_n) \geq 1/2$ and the second term in parentheses in \eqref{eq:Etilde} can be bounded as 
\begin{equation}\label{eq:Etilde1}
	\int_{x \in M} \lVert f_0 \rVert_2^2 \d\widetilde\pi(x) \leq \frac{V_n}{ \pi(T_n) \inf_{x \in T_n} \Vol B(x, \delta)} \leq  2V_n^{1+\eta_n},
\end{equation}
where $\eta_n  \to 0$ was defined at the beginning of the proof. Next,  consider the first term in parentheses in \eqref{eq:Etilde} for a given $\lambda_k \in \mathfrak S_C$. We argue as in the proof of Theorem \ref{th:D1}: using the Cauchy-Schwarz inequality and $\pi(T_n) \geq 1/2$ for $n$ large enough,
\begin{align*}
	\int_{x \in M} \langle f_0, \phi_k \rangle^2  \d\widetilde\pi(x) &\leq 2 \int_{x \in T_n}\left( \int_{y \in M} f_0^2(y) \d\pi(y) \int_{y\in M}\phi^2_k(y)\mathds{1}_{B(x, \delta)}(y) \d\pi(y) \right)\d \pi(x) \\&\leq 2 V_n \int_{x \in T_n} \int_{y \in M} \frac{\mathds 1_{B(x, \delta)}(y)}{\Vol B(x, \delta)} \phi_k^2(y) \d\pi(y)\d \pi(x)\\&\leq  2 V_n \int_{y \in M} \phi_k^2(y) \int_{x \in T_n} \frac{\mathds 1_{B(y, \delta)}(x)}{\Vol B(x, \delta)} \d \pi(x)\d\pi(y)\\&\leq  2 V_n^{\eta_n} e^{\delta(d-1) }.
\end{align*}

We recall that $t \geq (1+\varepsilon)t_\star$. By assumption, there exists $\gamma > 0$ such that $s_1 \leq 1- \gamma$ for all $n$ large enough. We then write
%, and from now on, $C$ is a positive constant which depends on $\delta, d$ but not on $n$, and which can change from one line to line. 
\begin{align*}
	\sum_{ \lambda_k \in \mathfrak S_C } e^{s_{k} (d-1) t} & =  \sum_{\lambda_k \in \mathfrak S_C} (d-1)  t \int_{-\infty}^{1-\gamma} \mathds{1}_{u \leq s_{k}} e^{u(d-1)t}\d u
	\\&\leq (d-1) t\int_{-\infty}^{1-\gamma} \# \left\lbrace \lambda \in \mathfrak S_C \vert \lambda \leq \sigma_d^2 (1-u^2)\right\rbrace e^{u(d-1)t}
	\d u.
\end{align*}
The Sarnak-Xue density property implies that for all $u\in [0,1)$,
$$
\# \left\lbrace \lambda \in \mathfrak S_C \vert \lambda \leq \sigma_d^2 (1-u^2)\right\rbrace \leq V_n^{1-u + \eta'_n}.
$$
where $\eta'_n$ is a vanishing sequence. We deduce that
\begin{align*}
	\sum_{ \lambda_k \in \mathfrak S_C } e^{s_{k} (d-1) t}  &\leq (d-1) t  \left( \int_{-\infty}^0  V_n^{1 + \eta'_n} e^{u(d-1)t} \d u  + \int_{0}^{1  - \gamma}  V_n^{1-u+\eta'_n} e^{u(d-1)t } \d u \right)
	\\& =  V_n^{1 + \eta'_n} + \frac{(d-1) t V_n^{1 + \eta'_n}}{(d-1) t - \ln V_n}  \left( V_n^{\gamma-1} e^{(1-\gamma)(d-1)t} - 1\right)\\
	& \leq \frac{V_n^{\gamma + \eta'_n}(1+\varepsilon)}{\varepsilon}  e^{(1-\gamma)(d-1)t}, 
\end{align*}		
where in the last line, we have used that $\ln V_n \leq (d-1) t/(1+\varepsilon)$. 
All in all, in \eqref{eq:Etilde}, combining our last bound with \eqref{eq:Etilde1}, we find for some new vanishing sequence $\eta''_n$ (depending on $\varepsilon,d$) that
\begin{align*}
	\E_{\widetilde\pi} \sup_{s \geq t} \DTV\left( \mu_s, \pi\right)^2 &\leq V_n^{1+\eta''_n} t e^{-(d-1)t}  + V_n^{\gamma+\eta''_n}  t e^{-\gamma (d-1)t} \\
	& \leq  C V_n^{\eta''_n -\gamma \varepsilon/2}, 
\end{align*}
where $C$ is such that $2 t \leq C \exp ( \gamma (d-1) t   \varepsilon /2)$ for all $t \geq 0$ and we have used $t \geq (1+\varepsilon) t_\star$. Finally, using Markov's inequality, we get
\[\widetilde\pi  \left(\left\{ x :   \sup_{s \geq t} \DTV(\mu_s,\pi) \geq V_n^{-\gamma \varepsilon / 8} \right\}  \right) 
\leq \frac{ \E_{\widetilde\pi} \sup_{s \geq t} \DTV \left( \mu_s, \pi\right)^2}{V_n^{-\gamma \varepsilon / 4} }
\leq   C V_n^{\eta''_n -\gamma \varepsilon /4}.
\]
The right side goes to $0$. We take $t = (1+\varepsilon)t_\star$. This concludes the proof of Theorem \ref{th:S2}. \qed

\section{Brownian Motion} \label{BM}

\subsection{Hyperbolic Brownian motion}

Let $(\widetilde W_t)_{t \in \Rp}$ be the Brownian motion on the hyperbolic space $\dH^d$ starting at $o = (0,\dots,0,1) \in \dH^d$ with infinitesimal generator $\Delta_{\dH^d} / 2$. 
Let $d (x,y)$ be the distance in $\dH^d$. It is a basic fact that $\widetilde W_t$ is radially symmetric. That is, the density of the law of $\widetilde W_t$ with respect to the hyperbolic measure $\ell_d$ at $x$ depends only on $d(x,o)$, see for example \cite[Subsection 3.3]{MR1882015}. One consequence is the following. Let $\widetilde D_t = d ( \widetilde W_t , o)$  and $\widetilde Z_t = (\widetilde X_t, \widetilde v_t) = \gamma_t (o , v_0) $, $t \in \Rp$, be the geodesic process started at $\widetilde X_0 = o$, independent of $\widetilde B$. We have the identity in law for any $t\in \Rp$:
\begin{equation*}\label{eq:B2X}
	\widetilde W_t \stackrel{d}{=} \widetilde X_{\widetilde D_t}.
\end{equation*} 
As explained in the proof of Proposition \ref{prop:Croissancesupport}, through the covering map the above equation remains true on any compact hyperbolic manifold $M$: for any $t\in \Rp$:
\begin{equation}\label{eq:B2X2}
	W_t \stackrel{d}{=}  X_{\widetilde D_t},
\end{equation}
where $(W_t)$ is the Brownian motion on $M$ started at $x$, $(X_t)$ is the position of the geodesic process on $M$ started at $x$ and $\widetilde D$ is as above and independent of $(X_t)_{t \in \Rp}$. We can also randomly sample the initial condition $x$: if $\mu_0$ is a probability distribution on $M$, then \eqref{eq:B2X2} holds where $(W_t)$ is the Brownian motion on $M$ with initial condition $\mu_0$ and $(X_t)$ is the position of the geodesic process on $M$ started with initial condition $\mu_0$.

We shall also use the fact that $\widetilde W_t$ has an explicit rate of escape: almost surely, 
\begin{equation}\label{eq:rate}
	\lim_{t \to \infty }  \frac{\widetilde D_t}{t} =  \frac{d -1}{2}, 
\end{equation}
see \cite{MR3176419,MR3631821}.

\subsection{Proof of Corollary \ref{coro:Brownianmixing1} and Corollary \ref{coro:Brownianmixing2}.}

We set $\beta = \frac{d-1}{2}$ and $t_\star = t_\star^n = \frac{\ln V_n}{d-1}$. We denote by $W_t = W_t^n$ the Brownian motion on $M_n$. We start with the upper bound. We fix $\eps > 0$ and set $s =  (1 + \eps) 2 \ln V_n = (1+ \eps) t_\star / \beta $. From the Markov property, it is standard that the map $t \to \DTV(\cL( W_{x,t}) , \pi)$ is non-increasing (see, for example, \cite[Proposition 3.3]{MR2375599}). Since $\eps >0$ is arbitrary, to prove the upper bound of Corollary \ref{coro:Brownianmixing1} and Corollary \ref{coro:Brownianmixing2}, it is thus sufficient to check that 
\begin{equation}\label{eq:ubB}
	\lim_{n \to \infty} \DTV( \cL( W_{x,s}) , \pi )  = 0,  
\end{equation}
where $x = x_n$ is as in Theorem \ref{th:S1} or as in Theorem \ref{th:S2} respectively. Fix some $\delta \in (0,1]$ as in the proof of Theorem \ref{th:S1} and Theorem \ref{th:S2} and let 
\begin{equation}\label{eq:choicerho}
	0 < \rho \leq \min ( \delta, \rho(x)).
\end{equation}

For any $\alpha \in (0,1)$, there exists $c_0 = c_0(\alpha) >0$ such that with $t_0 = c_0 \rho ^2$, we have 
\begin{equation}\label{eq:condBt}
	\dP ( \widetilde W_{t_0} \in B_{\dH^d} (o,\rho ) )\geq 1 - \alpha. 
\end{equation}
Indeed, according to \cite[Theorem 3.1]{Davies}, there exists $C >0$ which depends only on $d$ such that for all $t$ in $(0,1]$, for all $r \geq 0$, \begin{align}
	p_t(r) &\leq \frac{C}{t^{d/2}}e^{-\sigma_d t-r^2/4t-(d-1)r/2}(1+r+t)^{(d-3)/2}(1+r) \label{eq:bbrownien}
	\\& \leq \frac{C}{t^{d/2}}e^{-r^2/4t-(d-1)r/2}(1+r+t)^{(d-3)/2}(1+r).
\end{align}

For $r\geq 1$ and $t \leq 1$, $e^{-(d-1)r/2}(1+r+t)^{(d-3)/2}(1+r)$ is uniformly bounded from above by some $C'>0$ which depends only on $d$. So for all $0<t\leq 1$, for some constant $C$ which changes from one line to another, but only depends on $d$, \begin{align*}
	\p(\widetilde W_t \notin B_{\dH^d}(o, 1)) &\leq \frac{C}{t^{d/2}} \int_1^\infty e^{-r^2/4t}\sinh^{d-1} r\ \d r
	\\ & \leq \frac{C}{t^{d/2}} \int_1^\infty e^{-r^2/4t + (d-1)r}\ \d r
	\\ &\leq \frac{C}{t^{d/2}} \int_1^\infty e^{-(r/(2\sqrt{t}) - (d-1)\sqrt{t})^2}\ \d r
	\\   &\leq \frac{C}{t^{(d-1)/2}} \int_{1/(2\sqrt{t}) - (d-1)\sqrt{t}}^\infty e^{-x^2}\ \d x.
\end{align*}
As $t\to 0$, $\frac{1}{2\sqrt{t}} - (d-1)\sqrt{t} \to \infty$ and so \[\int_{1/(2\sqrt{t}) - (d-1)\sqrt{t}}^\infty e^{-x^2}\ \d x \sim \frac{e^{-\left(\frac{1}{2\sqrt{t}} - (d-1)\sqrt{t}\right)^2}}{\frac{1}{2\sqrt{t}} - (d-1)\sqrt{t}}.\]
In particular, $\p(\widetilde W_t \notin B_{\dH^d}(o, 1)) \limit{}{t\to 0} 0$ so for $t$ smaller than some $c(\alpha)$, \begin{equation}\label{eq:alpha1}
	\p(\widetilde W_t \notin B(o, 1)) \leq \alpha/2.
\end{equation}
In addition, from $\eqref{eq:bbrownien}$, we get that for any $0 < t, r \leq 1$, the radial density of $\widetilde W_t$ has a Radon-Nikodym derivative with respect to the radial density of the Euclidean Brownian motion $B_t$ in $\R^d$ -- that is $\frac{\Gamma_{d} r^{d-1}}{t^{d/2}}e^{-\frac{r^2}{4t}}$ -- which is uniformly bounded from above by some $C>0$ which depends only of $d$. In particular, \[\p(\widetilde W_t \in B_{\dH^d}(o, 1)\setminus B_{\dH^d}(o, \rho)) \leq C\p(B_t \in B_{\R^d}(0, 1)\setminus B_{\R^d}(0, \rho)).\]
As basic consequence of scale invariance of the Brownian on $\R^d$ -- $W_t \stackrel{d} {=} \sqrt t W_1$ -- we can find $c'(\alpha)$ such that for any $t \leq c'(\alpha)\rho^2$, this is smaller than $\alpha/2$. Combined with $\eqref{eq:alpha1}$, it shows that it suffices to take $t_0 = \min(c, c')\rho^2$. Let now $\widetilde g_{0}$ be the density of $\widetilde W_{t_0}$ with respect to the hyperbolic measure $\ell_d$. Again, because $\rho \leq 1$, by $\eqref{eq:bbrownien}$, we have 
$$
\sup_{x \in B_{\dH^d} (o,\rho)} \widetilde g_{0} (x) \leq c_1  t_0^{-d/2} =  c_2   \rho ^{-d},
$$
where $c_1 = c_1(\alpha)$ is some constant and $c_2 = c_1 c_0^{-d/2}$.

Note than in $M_n$, $W_{x,t_0}$ can be obtained by projecting $\widetilde W_{t_0}$ (where $o$ is projected on $x$). We keep this point of view and consider $\hat g_0$ the density of $W_{x,t_0}$ conditioned on the event $$E = \{\widetilde W_{t_0}  \in B_{\dH^d} (o, \rho) \}$$ with respect to the uniform measure on $B(x,\delta)$: $\IND ( y \in B(x,\delta) ) / \Vol B(x,\delta) d \Vol(y)$. From what precedes and \eqref{eq:condBt}, we find that
\begin{equation}\label{eq:condBt2}
	\sup_{y \in B(x,\rho) } \hat g_0(x) \leq \frac{c_2 \Vol( B(x,\delta) )}{1 -\alpha} \rho^{-d} \leq c \rho^{-d},
\end{equation}
with $c = c(\alpha) = c_2 \Vol_{\mathds H^d} ( B(o,\delta) )  / (1-\alpha)$.

We are ready for the proof of the upper bound. We may find $b < \beta$  such that $s b \geq (1 + \eps / 2) t_\star$. For any measurable event $A$, from \eqref{eq:condBt}, we have 
\begin{equation}\label{eq:bdAA}
	\p ( W_{x,s + t_0} \in A ) - \pi(A) \leq  \p(E) \left( \p ( W_{x,s+t_0} \in A | E )   - \pi(A) \right)   + 2 \alpha,
\end{equation}
where $\p( A | B) = p(A \cap B) / \p(B)$ denotes the conditional probability.   From the Markov property and \eqref{eq:condBt2}, we find
\begin{align*}
	\p ( W_{x,s+t_0} \in A | E )   - \pi(A) & =  \int_{M_n}  \left( \dP ( W_{y,s} \in A ) - \pi(A) \right) \hat g_0(y)  f_0(y) d \Vol (y) \\
	&\leq  c \rho^{-d}  \int_{M_n}   \left| \dP ( W_{y,s} \in A ) - \pi(A) \right| f_0(y) d \Vol (y) \\ 
	& = c \rho^{-d}     \left| \dP ( W_{X_{x,0},s} \in A ) - \pi(A) \right|,
\end{align*}
where $(W_{X_{x,0},t})_{t \geq 0}$ is the Brownian motion on $M_n$ with initial condition $X_{x,0}$ uniform on $B(x,\delta)$.  As in Theorem \ref{th:S1} or Theorem \ref{th:S2}, let $(X_{x,t})$ be the geodesic process on $M_n$ with initial condition uniform on $B(x,\delta)$. From \eqref{eq:B2X2}, if $\nu_t$ is the law of $\widetilde D_t$, we have, 
\begin{align*}
	\p ( W_{X_{x,0},s} \in A ) - \pi(A) & = \int_0 ^\infty \p ( X_{x,t} \in A) \nu_{s} (\d t) - \pi(A) \\
	&  \leq 2 \p (  \widetilde D_{s} <  s b ) + \int_{( 1+ \eps / 2 ) t_\star}^\infty (\p(X_{x,t} \in A) - \pi(A) ) \nu_{s} (\d t) \\
	&   \leq 2 \p (  \widetilde D_{s} <  s b ) + \sup_{t \geq ( 1+ \eps / 2 ) t_\star} \DTV( \cL (X_{x,t}),\pi). 
\end{align*}
If  $ x  = x_n$ is as in Theorem \ref{th:S1} or Theorem \ref{th:S2} respectively, we can take in \eqref{eq:choicerho} $\rho = V^{-\eta_n}_n$ for some sequence $\eta_n \to 0$. The second term of the last shown inequality goes to $0$ while the first term goes to $0$ from \eqref{eq:rate}. Moreover, these terms go to $0$ at a rate faster than $\rho^{-d} = V_n^{d \eta_n}$ as can be seen from the proofs of Theorem \ref{th:S1} or Theorem \ref{th:S2} for the second term and from \cite[Proposition 3.2]{MR3176419}  for the first term. Taking the supremum over all sets $A$, from \eqref{eq:bdAA}, we thus have checked that 
$$
\limsup_{n \to \infty} \DTV( \cL (W_{x,s + t_0}),\pi) \leq 2 \alpha. 
$$ 
Since $t_0 = O(1)$ and $\alpha$ can be taken arbitrarily small, it concludes the proof of the upper bound \eqref{eq:ubB}.

It remains to prove the corresponding lower bound. We argue as in the proof of Proposition \ref{prop:LB}. We fix $0 < \eps < 1$ and set $s' =  (1 - \eps) 2 \ln V_n = (1 -  \eps) t_\star / \beta $. There exists $b'>\beta$ such that $s' b' \leq (1 - \eps / 2) t_\star$. Since  $t \to \DTV(\cL( W_{x,t}) , \pi)$ is non-increasing, it suffices to prove that 
\begin{equation}\label{eq:lbB}
	\lim_{n \to \infty} \sup_{x \in M_n} \DTV( \cL( W_{x,s'}) , \pi )  = 1.  
\end{equation}

{Take $x \in M_n$. For $r > 0$, we denote by $B(x,r)$ the closed ball in $M_n$ with radius $r$ and center $x$. Similarly, we denote by $\widetilde B(r)$ the closed ball in $\dH^d$ with center $o$ and radius $r$. Arguing as in the proof of Proposition \ref{prop:LB}, by using the covering map, we have, for some $C >0$, 
	\begin{equation}\label{eq:upvolball}
		\Vol( B(x,r) ) \leq \ell_d ( \widetilde B (r) ) \leq  C \exp ( (d-1) r ).  
	\end{equation}
	Furthermore, the distance in $M_n$ between $x$ and $W_{x,t}$ is stochastically dominated by $\widetilde D_t = d(\widetilde W_{t},o)$. In other words, for any $r,t$,
	$$
	\dP ( W_{x,t} \in B ( x, r  ) ) \geq \dP ( \widetilde D_t  \leq  r  ).
	$$
	In particular,  
	$$
	\dP ( W_{x,s'} \in B ( x, (1-\eps/2)t_\star  ) ) \geq \dP ( \widetilde D_{s'}  \leq  b' s' ).
	$$  
	On the other hand, from \eqref{eq:upvolball}, we have 
	$$
	\pi( B ( x, (1-\eps/2)t^n_\star  ) =  \frac{\Vol (  B ( x, (1-\eps/2)t_\star )}{V_n} \leq C V_n^{-\eps /2}. 
	$$
	From the definition of $\DTV$, it follows that 
	\begin{align*}
	\DTV( \cL( W_{x,s'}) , \pi  ) &\geq \dP ( W_{x,s'} \in B ( x, (1-\eps/2)t_\star  ) ) - \pi( B ( x, (1-\eps/2)t_\star  ) \\&\geq \dP ( \widetilde D_{s'}  \leq  b' s' )  - C V_n^{-\eps /2}.
	\end{align*}
	Using \eqref{eq:rate} and $b' > \beta$, it concludes the proof of \eqref{eq:lbB}.
}
\qed

\paragraph{Formal statements}
The authors have no conflict of interest to declare. There is no data associated to this manuscript.

\printbibliography

\end{document}